\documentclass[12pt]{article}

\usepackage{latexsym}
\usepackage{amsfonts}
\usepackage{amssymb}
\usepackage{amsmath}
\usepackage{mathrsfs}
\usepackage{hyperref}

\newcommand{\be}{\begin{equation}}
\newcommand{\ee}{\end{equation}}
\newcommand{\bea}{\begin{eqnarray}}
\newcommand{\eea}{\end{eqnarray}}
\newcommand{\bean}{\begin{eqnarray*}}
\newcommand{\eean}{\end{eqnarray*}}
\newcommand{\brray}{\begin{array}}
\newcommand{\erray}{\end{array}}

\newtheorem{dfn}{Definition}[section]
\newtheorem{thm}[dfn]{Theorem}
\newtheorem{lmma}[dfn]{Lemma}
\newtheorem{ppsn}[dfn]{Proposition}
\newtheorem{crlre}[dfn]{Corollary}
\newtheorem{xmpl}[dfn]{Example}
\newtheorem{rmrk}[dfn]{Remark}

\newcommand{\bdfn}{\begin{dfn}\rm}
\newcommand{\bthm}{\begin{thm}}
\newcommand{\blmma}{\begin{lmma}}
\newcommand{\bppsn}{\begin{ppsn}}
\newcommand{\bcrlre}{\begin{crlre}}
\newcommand{\bxmpl}{\begin{xmpl}}
\newcommand{\brmrk}{\begin{rmrk}\rm}

\newcommand{\edfn}{\end{dfn}}
\newcommand{\ethm}{\end{thm}}
\newcommand{\elmma}{\end{lmma}}
\newcommand{\eppsn}{\end{ppsn}}
\newcommand{\ecrlre}{\end{crlre}}
\newcommand{\exmpl}{\end{xmpl}}
\newcommand{\ermrk}{\end{rmrk}}

\newcommand{\cla}{\mathcal{A}}
\newcommand{\clb}{\mathcal{B}}

\newcommand{\clh}{\mathcal{H}}

\newcommand{\clg}{\mathcal{G}}

\author{S. Sundar}
\title{Toeplitz algebras associated to Endomorphisms of Ore semigroups}

\begin{document}
\maketitle 
\begin{abstract}
In this paper, we consider the Toeplitz algebra associated to actions of Ore semigroups on $C^{*}$-algebras. In particular, we consider injective and surjective actions of such semigroups. We use the theory of groupoid dynamical systems to represent the Toeplitz algebra as a groupoid crossed product. We also discuss the K-theory of the Toeplitz algebra in some examples. For instance, we show that for the semigroup of positive matrices, the K-theory of the associated Toeplitz algebra vanishes.
\end{abstract}

\noindent {\bf AMS Classification No. :} {Primary 22A22; Secondary 54H20, 43A65, 46L55.}  
\tableofcontents

\section{Introduction}
\footnote{\textbf{Acknowledgements}: 
This work should be seen as  a continuation of my work done in colloboration with Prof. Renault. I would like to thank Prof. Renault for inviting me for a month at Orleans during Sept-Oct 2013. For, the results of this paper are shaped by the discussions that I had with Prof. Renault on groupoids. I would also like to thank Prof. Emmanuel Germain of university of Caen for funding my visit and Anne Liger for taking care of everything at Orleans.}

Group $C^{*}$-algebras and the associated crossed products have always been at the centre of interest in the theory of $C^{*}$-algebras. In the nineties, Gerard J. Murphy in a series of papers studied $C^{*}$-algebras associated to semigroups and semigroup actions. The reader is referred to \cite{Murphy96}, \cite{Murphy94} and the references therein for the work of Murphy on semigroup $C^{*}$-algebras. The study of semigroup $C^{*}$-algebras is revived again in the last few years with the work(\cite{Cuntz}) of Cuntz on the $C^{*}$-algebras associated to the "ax+b"-semigroup.  We  refer to \cite{Li-semigroup}, \cite{Li13} for recent developments on semigroup $C^{*}$-algebras. Even though we can now say that the literature on semigroup $C^{*}$-algebras and semigroup crossed products is vast,  the literature on $C^{*}$-algebras associated to  topological semigroups is few in comparision to that on discrete semigroups. 

The paper \cite{Renault_Muhly} by Muhly and Renault on Wiener-Hopf $C^{*}$-algebras associated to cones in $\mathbb{R}^{n}$ can be considered as the first systematic attempt in the direction of topological semigroups.  
Renault and Muhly made extensive use of groupoid techniques to explore the Wiener-Hopf algebras associated to polyhedral cones and self-dual cones. Later the works of Nica (\cite{Nica_WienerHopf}, \cite{Nica90}), Hilgert and Neeb (\cite{Hilgert_Neeb}, \cite{Hilgert_Neeb1}) extended the analysis initiated in \cite{Renault_Muhly} to other semigroups in particular, by Hilgert and Neeb, to subsemigroups of Lie groups. In \cite{Jean_Sundar}, based on the techniques used in \cite{Hilgert_Neeb}, the Wiener-Hopf $C^{*}$-algebras associated to  Ore semigroups are studied and a groupoid picture is obtained for them. Nevertheless, up to the author's knowledge, the only reference in the literature where Wiener-Hopf or Toeplitz type $C^{*}$-algebra associated to  actions of topological semigroups is considered is \cite{KS97}. The action of the semigroup of the positive real line on a $C^{*}$-algebra is considered by  Khoshkam and Skandalis in \cite{KS97}. In fact, the present paper is an outgrowth of the author's effort to understand \cite{KS97}, the author's desire to place \cite{KS97} within the framework of groupoids and to extend some of the results in \cite{KS97} to possibly more general semigroups. This article owes its existence to \cite{KS97}. In fact, the title of this paper is inspired by that of  \cite{KS97}.

A brief outline of the problem considered in this paper is given below. First let us recall the construction of the reduced crossed product. Let $G$ be a locally compact group, $A$ be a $C^{*}$-algebra and $\alpha:G \to Aut(A)$ be a strongly continuous action. Consider the Hilbert $A$-module $E:=A \otimes L^{2}(G)=L^{2}(G,A)$. Let $\mathcal{L}_{A}(E)$ be the $C^{*}$-algebra of adjointable operators on $E$. Define $\widetilde{\pi}:A \to \mathcal{L}_{A}(E)$ and $\widetilde{U}:G \to \mathcal{L}_{A}(E)$ by the formula:
\begin{align*}
\widetilde{\pi}(x)(\xi)(h):&=\alpha_{h}^{-1}(x)\xi(h) \\
\widetilde{U}_{g}(\xi)(h):&= \xi(g^{-1}h)
\end{align*}
for $x \in A$, $g \in G$ and $\xi \in L^{2}(G,A)$. The reduced crossed product $A \rtimes_{red} G$ is defined as the $C^{*}$-algebra generated by $\{\int \widetilde{\pi}(f(g))\widetilde{U}_{g}dg: f \in C_{c}(G,A)\}$. Here $C_{c}(G,A)$ denotes the algebra of compactly supported continuous functions  on $G$ taking values in $A$.  The formula for $\widetilde{\pi}$ makes use of the fact that $\alpha_{g}$ is invertible for every $g \in G$.  Thus, it is not clear how to define the reduced crossed product $A \rtimes_{red} G$ when $G$ is just a semigroup and the action is only by endomorphisms. However the representations $\widetilde{\pi}$ and $\widetilde{U}$ can be twisted using the unitary on $L^{2}(G)$ induced by the inversion operation $G \ni g \to g^{-1}\in G$ up to the modular function. Thus define $\pi:A \to \mathcal{L}_{A}(E)$ and $U:G \to \mathcal{L}_{A}(E)$ by the formulas:
\begin{align*}
\pi(x)(\xi)(h):&=\alpha_{h}(x)\xi(h) \\
U_{g}(\xi)(h):&= \xi(hg^{-1})\Delta(g)^{-\frac{1}{2}} 
\end{align*}
Let $V$ be the unitary on $L^{2}(G,A)$ defined by $V\xi(g):=\xi(g^{-1})\Delta(g)^{-\frac{1}{2}}$. Then $V\widetilde{\pi}(.)V^{*}=\pi(.)$ and $V\widetilde{U_{g}}V^{*}=U_{g^{-1}}$ for $g \in G$. Thus $A \rtimes_{red} G$ can be defined as the $C^{*}$-algebra generated by $\{\int \pi(f(g))U_{g}dg: f \in C_{c}(G,A)\}$. The representation $\pi$ and $U$ makes sense as long as $G$ is just a subsemigroup of a locally compact group and  $\alpha$ is just an action by endomorphisms.

Let $P \subset G$ be a closed subsemigroup of a locally compact group $G$ containing the identity element. In this paper, we consider only semigroups of Ore type. We assume that $P$ is right Ore i.e.  $PP^{-1}=G$ and $P$ is solid i.e. $Int(P)$ is dense in $P$. At present, it is not clear to the author how to extend the analysis of this paper to subsemigroups which are not of Ore type. Nevertheless, this covers quite a few examples like polyhedral cones in $\mathbb{R}^{n}$, cones of positive matrices in $M_{n}(\mathbb{C})$, the subsemigroup of positive elements in the Heisenberg group etc. Let $P \subset G$ be a closed right Ore , solid subsemigroup of $G$ . Let $\alpha:P \to End(A)$ be a strongly continuous homomorphism where $End(A)$ denote the semigroup of $*$-algebra endomorphisms of $A$. 
Consider the Hilbert $A$-module $E:=A \otimes L^{2}(P)=L^{2}(P,A)$. Let $\pi:A \to \mathcal{L}_{A}(E)$ and $V:P \to \mathcal{L}_{A}(E)$ be defined by the formulas:
\begin{align*}
\pi(x)(\xi)(a):&=\alpha_{a}(x)\xi(a) \\
V_{b}(\xi)(a):&= \xi(ab^{-1})\Delta(b)^{-\frac{1}{2}}.
\end{align*}    
Then $\pi$ is a representation and $V$ is a strictly continuous representation of $P$ as isometries. Moreover we have the covariance relation
\[
V_{a}^{*}\pi(x)V_{a}=\pi(\alpha_{a}(x))
\] 
for $a \in P$ and $x \in A$. Just like in the group case, we define  the $C^{*}$-algebra generated by $\{\int \pi(f(a))V_{a} da: f \in C_{c}(P,A)\}$ as the reduced crossed product of $A$ by $P$ and denote it by $A \rtimes_{red,\alpha} P$.
There is another related $C^{*}$-algebra which we call the Wiener-Hopf or the Toeplitz algebra associated to the action $\alpha$. For $g \in G$, write $g=ab^{-1}$ with $a,b \in P$ and let $W_{g}:=V_{b}^{*}V_{a}$. Then it is trivial to verify that $W_{g}$ is well defined and the map $G \ni g \to W_{g}\in \mathcal{L}_{A}(E)$ is strictly continuous. We denote the $C^{*}$-algebra generated by $\{\pi(f(g))W_{g} dg: f \in C_{c}(G,A)\}$ as $\mathfrak{T}(A,P,\alpha)$. Note that if $A=\mathbb{C}$, then $\mathfrak{T}(A,P,\alpha)$ is just the usual Wiener-Hopf algebra whose study was initiated in \cite{Renault_Muhly}. The Toeplitz algebra $\mathfrak{T}(A,P,\alpha)$ is studied in \cite{KS97} when $P$ is either $\mathbb{N}$ or the additive semigroup $[0,\infty)$ of positive reals.

In this paper, we examine the $C^{*}$-algebra $A \rtimes_{red,\alpha} P$ in the  following two cases
\begin{enumerate}
\item[(1)] The algebra $A$ is unital and the action $\alpha$ is unital and injective i.e. for each $a \in P$, $\alpha_{a}$ is unital and injective.
\item[(2)] The action $\alpha$ is surjective i.e. $\alpha_{a}$ is surjective for each $a \in P$.
\end{enumerate}
We show that $A \rtimes_{red,\alpha} P$ coincides with $\mathfrak{T}(A,P,\alpha)$ in both the injective and the surjective case. The main result(s) of this paper is to represent $\mathfrak{T}(A,P,\alpha)$ as a groupoid crossed product in both the injective and the surjective case (Thm. \ref{main_theorem} and Thm. \ref{surjective main_theorem}).

Let $\mathcal{C}(G)$ be the space of closed subsets of $G$ endowed with the Fell topology. With respect to the Fell topology, $\mathcal{C}(G)$ is compact. Moreover $G$ acts on $\mathcal{C}(G)$ on the right by right multiplication. Denote the closure of $\{P^{-1}a: a \in P\}$ by $\Omega$. Then $\Omega$ is compact and $P$-invariant. Let $\clg$ be the groupoid obtained by restricting the transformation groupoid $\mathcal{C}(G) \rtimes G$ to $\Omega$. The groupoid $\clg$ is the same as the Renault-Deaconu semi-direct product groupoid $\Omega \rtimes P$ which we also call as the Wiener-Hopf groupoid. It is shown in \cite{Jean_Sundar} that the Wiener-Hopf groupoid posseses a Haar system, one such is given by the restriction of the Haar measure on $G$  to the fibres. 

In both the injective and the surjective case, we construct an upper semi-continuous bundle $\cla$ over the unit space $\Omega$ and a left action of the Wiener-Hopf groupoid $\clg:=\Omega \rtimes P$ on $\cla$. 
However the construction of the bundle $\cla$ is different in both  cases and at present the author is not able to uniformise the construction. The construction of the bundle $\cla$ (its fibres) in the injective case relies on a dilation principle whereas in the surjective case it relies on a quotient construction. But in both  cases, we show that $\mathfrak{T}(A,P,\alpha)$ is isomorphic to the reduced crossed product $\cla \rtimes_{red} \clg$. This along with the fact that $\Omega \rtimes P$ is equivalent to a transformation groupoid of the form $\widetilde{\Omega} \rtimes G$ implies that $A \rtimes_{red,\alpha} P$ is Morita-equivalent to a reduced crossed product of the form $B \rtimes_{red} G$ for some $C^{*}$-algebra. This Morita equivalence is then used to derive some $K$-theoretic consequences. For example, we prove that if $P$ is the cone of positive matrices then $K_{i}(A \rtimes_{red} P)=0$. This generalises  the classic result that if $\alpha:\mathbb{R} \to Aut(A)$ is a strongly continuous action then the Toeplitz extension of $A \rtimes_{\alpha} \mathbb{R}$ has trivial $K$-theory. The vanishing of the $K$-theory of the Toeplitz extension of $A \rtimes_{\alpha} \mathbb{R}$ forms the basis of Rieffel's proof of Connes-Thom isomorphism.

This paper is organised as follows.

In Section 2, we recall the necessary preliminaries regarding groupoid dynamical systems. Even though excellent references are availabe regarding groupoid dynamical systems, \cite{KS04} and \cite{MW08} to mention two, I believe that devoting two pages to recall the basics will make this paper more readable.

In Section 3, we recall the construction of the Wiener-Hopf groupoid and a few results  from \cite{Jean_Sundar}. In Section 4, we formally define the Toeplitz algebra and the reduced crossed product associated to an action of a right Ore semigroup on a $C^{*}$-algebra.

Sections 5 and 6 form the heart of this paper. In Section 5, we further make the assumption that $P$ is left Ore i.e. $P^{-1}P=G$. We construct, in Section 5, a groupoid dynamical system $(\cla,\clg,\alpha)$ associated to a unital injective action $\alpha:P \to End(A)$. We make use of the dilation principle which roughly says that injective actions of Ore semigroups can be dilated to bijective actions. 

In Section 6, we treat the surjective case where the groupoid dynamical system $(\cla,\clg,\alpha)$ is obtained using a quotient construction. We show in Sections 5 and 6 that the Toeplitz algebra $\mathfrak{T}(A,P,\alpha)$ is isomorphic to the reduced crossed product $\cla \rtimes_{red} \clg$. We also show that $\mathfrak{T}(A,P,\alpha)$ coincides with $A \rtimes_{red} P$.

We use the results of Sections 5 and 6 in Section 7 to prove that in both the surjective and the injective case, the reduced crossed product $A \rtimes_{red} P$ is Morita-equivalent to a genuine crossed product of the form $B \rtimes_{red} G$ for some $C^{*}$-algebra $B$. Using this Morita equivalence, we derive some $K$-theoretic consequences in Section 8 under the assumption that $G$ satisfies the Connes-Thom isomorphism and $P$ is connected. We must mention that the Morita equivalence derived in Section 7 can be used to obtain $K$-theoretic results for discrete semigroups like $\mathbb{N}$. But we don't pursue this case in this paper.

We conclude, in Section 9, by discussing the example of the semigroup of positive elements in a special Euclidean Jordan algebra i.e. real vector subspaces of self-adjoint matrices which contains the identity and closed under taking the anti-commutator of two elements. We make use of the Cayley transform to embedd self-adjoint matrices in the unitary group and extend the action of self-adjoint matrices on itself by addition to an action  on the unitary group. We obtain an explicit picture of the Wiener-Hopf groupoid in this case and use it to prove that $K_{i}(A \rtimes_{red} P)=0$ where $P$ is the semigroup of positive elements in a  special Euclidean Jordan algebra.

\section{Groupoid dynamical systems}
In this section, we recall the necessary preliminaries regarding groupoid dynamical systems. For more details and proofs, we refer to \cite{KS04} and \cite{MW08}. We use the  notion  of an upper semi-continuous bundle of $C^{*}$-algebras over a topological space for which the reader is refered to \cite{Dana-Williams}. We also make use of vector valued integration for which \cite{Dana-Williams} is the reference. The author claims no originality of what follows in this section.

Let $X$ be a locally compact, Hausdorff topological space and let $\mathcal{A}$ be an upper semi-continuous bundle of $C^{*}$-algebras over $X$ with the projection map $p:\cla \to X$. For $x \in X$, let $A_{x}:=p^{-1}(x)$. The space of continuous sections of $\cla$ vanishing at infinity is denoted by $\Gamma_{0}(X,\cla)$.
For $f \in \Gamma_{0}(X,\cla)$, let\[ ||f||:=\sup_{x \in X}||f(x)||.\] Then $\Gamma_{0}(X,\cla)$ is a $C^{*}$-algebra where the algebraic operations are defined using the corresponding ones on the fibres $\{A_{x}\}_{x \in X}$. Moreover for $\phi \in C_{0}(X)$ and $f \in \Gamma_{0}(X,\cla)$, define $\phi f \in \Gamma_{0}(X,\cla)$ by the equation
\[
(\phi f)(x):=\phi(x)f(x) \in A_{x}.
\]  
Thus $\Gamma_{0}(X,\cla)$ is a $C_{0}(X)$ module. We denote the space of compactly supported continuous sections by $\Gamma_{c}(X,\cla)$. Clearly $\Gamma_{c}(X,\cla)$ is dense in $\Gamma_{0}(X,\cla)$. If $X$ is compact, we simply denote $\Gamma_{0}(X,\cla)$ by $\Gamma(X,\cla)$.
\begin{rmrk}
We will assume that all the upper semi-continuous bundles that we consider have enough sections i.e. for every $x \in X$, the map $\Gamma_{0}(X,\mathcal{A}) \ni f  \to f(x) \in A_{x}$ is onto. Moreover we will assume that $\cla$ is first countable and $X$ is locally compact, second countable and Hausdorff. Thus it is enough to consider sequences instead of nets. 
\end{rmrk}

Next we recall the notion of a groupoid dynamical system and the associated reduced crossed product. In what follows, $\mathcal{G}$ denotes a second countable, locally compact Hausdorff groupoid with a Haar system $\{\lambda^{x}: x \in \mathcal{G}^{(0)}\}$ where $\mathcal{G}^{(0)}$ is the unit space of $\mathcal{G}$. The range and source maps are denoted, as usual, by $r$ and $s$. For $x \in \clg^{(0)}$, we denote $r^{-1}(x)$ and $s^{-1}(x)$ by $\mathcal{G}^{(x)}$ and $\clg_{(x)}$ respectively.

\begin{dfn}
Let $\mathcal{A}$ be an upper semi-continuous bundle of $C^{*}$-algebras over $\mathcal{G}^{(0)}$. We say $\mathcal{G}$ acts on $\cla$ if for every $\gamma \in \mathcal{G}$, there exists a map $\alpha_{\gamma}:A_{s(\gamma)} \to A_{r(\gamma)}$ such that 
\begin{enumerate}
\item for $\gamma \in \mathcal{G}$, the map $\alpha_{\gamma}:A_{s(\gamma)} \to A_{r(\gamma)}$ is a $*$-algebra isomorphism,
\item if $r(\eta)=s(\gamma)$, $\alpha_{\gamma}\alpha_{\eta}=\alpha_{\gamma \eta}$  , and
\item for $\gamma \in \mathcal{G}$, $\alpha_{\gamma^{-1}}=\alpha_{\gamma}^{-1}$.
\item the map $\clg*\cla \ni (\gamma,a) \to \alpha_{\gamma}(a) \in \cla$ is continuous where  $\clg*\cla \subset \clg \times \cla$ is defined as $\clg*\cla:=\{(\gamma,a):a \in A_{s(\gamma)}\}$. The subset $\clg*\cla$ is given the subspace topology inherited from the product topology on $\clg \times \cla$.
\end{enumerate}  
We call the triple $(\mathcal{A},\mathcal{G}, \alpha:=\{\alpha_{\gamma}\})$ a groupoid dynamical system.
\end{dfn}

Let $(\mathcal{A},\mathcal{G},\alpha)$ be a groupoid dynamical system. Let  \[r^{*}\cla:=\{(\gamma,a ) \in \mathcal{G} \times \cla: a \in A_{r(\gamma)}\}.\] Then $r^{*}\cla$ is an upper semi-continuous bundle over $\mathcal{G}$ and is called the pull-back bundle corresponding to the range map $r:\clg \to \clg^{(0)}$.  We denote its space of continuous compactly supported sections by $\Gamma_{c}(\mathcal{G},r^{*}\cla)$. Thus 
\[
\Gamma_{c}(\mathcal{G},r^{*}\cla)=\{ f:\mathcal{G} \to \mathcal{A} : \textrm{ $f$ is continuous and $f(\gamma) \in A_{r(\gamma)}$ for $\gamma \in \mathcal{G}$}\}.
\]  For $f,g \in \Gamma_{c}(\mathcal{G},r^{*}\cla)$, let 
\begin{align*}
(f*g)(\gamma):&=\int f(\eta)\alpha_{\eta}(g(\eta^{-1} \gamma)) d\lambda^{r(\gamma)}, \\
 f^{*}(\gamma):&=\alpha_{\gamma}(f(\gamma^{-1})^{*}).
\end{align*}
With the multiplication and involution defined above, $\Gamma_{c}(\mathcal{G},r^{*}\cla)$ forms a $*$-algebra. Also $\Gamma_{c}(\mathcal{G},r^{*}\cla)$ forms a pre-Hilbert right $\Gamma_{c}(\mathcal{G}^{(0)},\cla)$ module where the right multiplication and the inner product are given as follows :

For $f,g \in \Gamma_{c}(\mathcal{G},\cla)$ and $\phi \in \Gamma_{c}(\mathcal{G}^{(0)},\cla)$, let 
\begin{align*}
(f.\phi)(\gamma):&=f(\gamma)\alpha_{\gamma}(\phi(s(\gamma))) \textrm{ for $\gamma \in \mathcal{G}$, and} \\
 \langle f,g \rangle(x):& = (f^{*}*g)(x) = \int \alpha_{\gamma}(f(\gamma^{-1})^{*})\alpha_{\gamma}(g(\gamma^{-1}))d\lambda^{(x)}(\gamma) \textrm{ for $x \in \mathcal{G}^{(0)}$}. 
\end{align*}

On completion, we obtain a  Hilbert $B:=\Gamma_{0}(\mathcal{G}^{(0)},\cla)$-module which we denote by $E$. The Hilbert module $E$ contains $\Gamma_{c}(\mathcal{G},r^{*}\cla)$ as  a dense subspace. For $f \in \Gamma_{c}(\clg,r^{*}\cla)$, let $\Lambda(f)$ be the operator on $\Gamma_{c}(\clg,r^{*}\cla)$ given by  $\Lambda(f)g=f*g$  for $g \in \Gamma_{c}(\clg,r^{*}\cla)$. One checks that for $f \in \Gamma_{c}(\clg,r^{*}\cla)$, $\Lambda(f)$ extends to an adjointable operator on $E$ which we still denote by $\Lambda(f)$. For a proof of this fact, see Section 3.6 of \cite{KS04}. In fact, $||\Lambda(f)|| \leq ||f||_{I}$ where 
\[
||f||_{I}:= \max \Big \{ \sup_{x \in \clg^{(0)}}\int ||f(\gamma)||d\lambda^{(x)}(\gamma), \sup_{x \in \clg^{(0)}} \int ||f(\gamma^{-1})|| d\lambda^{(x)}(\gamma) \Big \}.
\]

The map $\Lambda:\Gamma_{c}(\clg,r^{*}\cla) \to \mathcal{L}_{B}(E)$ is a faithful $*$-representation.  The reduced crossed product denoted $\mathcal{A} \rtimes_{red} \clg$ is defined as the closure of $\{\Lambda(f): f \in \Gamma_{c}(\mathcal{G},r^{*}\cla)\}$ in $\mathcal{L}_{B}(E)$. For $f \in \Gamma_{c}(\clg,r^{*}\cla)$, we let $||f||_{red}:=||\Lambda(f)||$. We also remark that the $I$-norm $||~||_{I}$ is a Banach algebra norm on $\Gamma_{c}(\clg,r^{*}\cla)$ and the involution is isometric.

\begin{rmrk}
For $f \in C_{c}(\clg)$ and $g \in \Gamma_{c}(\clg,r^{*}\cla)$, the convolution product $f*g$ defined by 
\[
f*g(\gamma):= \int f(\eta)\alpha_{\eta}(g(\eta^{-1}\gamma))d\lambda^{(r(\gamma))} 
\]
is a continuous compactly supported section of the bundle $r^{*}\cla$. Morover $||f*g||_{I} \leq ||f||_{I}||g||_{I}$ for $f \in C_{c}(\clg)$ and $g \in \Gamma_{c}(\clg,r^{*}\cla)$.
\end{rmrk}

For $x \in \mathcal{G}^{(0)}$, Let $\epsilon_{x}$ be the map $B:=\Gamma_{0}(\mathcal{G}^{(0)},\cla) \ni f  \to f(x) \in A_{x}$. For every $x \in X$, by Rieffel's Induction, one obtains a representation $\Lambda_{x}$ of $\cla \rtimes_{red} \clg$ on $E_{x}:= E \otimes_{B} A_{x}$ i.e $\Lambda_{x}(f)=\Lambda(f) \otimes 1$. 

Fix $x \in \clg^{(0)}$. Denote the external tensor product $L^{2}(\mathcal{G}^{(x)}) \otimes A_{x}$ by $L^{2}(\clg^{(x)},A_{x})$. The Hilbert $A_{x}$-module $L^{2}(\mathcal{G}^{(x)},A_{x})$ can also be regarded as the completion of $C_{c}(\clg^{(x)},A_{x})$ where for $\xi,\eta \in C_{c}(\clg^{(x)},A_{x})$, the inner product $<\xi,\eta>$ is given by \[
<\xi,\eta>(x):= \int \xi(\gamma)^{*}\eta(\gamma)d \lambda^{(x)}(\gamma)
\] 
for $x \in \clg^{(0)}$.
The map $\Gamma_{c}(\clg,r^{*}\cla)\otimes_{B} A_{x} \ni \xi \otimes a \to \widehat{\xi} \in C_{c}(\mathcal{G}^{(x)},A_{x})$ where $\widehat{\xi}(\gamma)=\alpha_{\gamma}(\xi(\gamma^{-1}))a$ extends to a unitary isomorphism. Thus we identify $E_{x}$ with $L^{2}(\mathcal{G}^{(x)},A_{x})$. Under this isomorphism, for $f \in \Gamma_{c}(\clg,r^{*}\cla)$ and $\xi \in C_{c}(\mathcal{G}^{(x)},A_{x})$, $\Lambda_{x}(f)\xi$ is given by
\[
(\Lambda_{x}(f)\xi)(\gamma):=\int \alpha_{\gamma}(f(\gamma^{-1}\gamma_{1}))\xi(\gamma_{1})d\lambda^{(x)}(\gamma_{1}).
\] 
The following are easily verifiable.

\begin{enumerate}
\item[(1)] For $f \in \Gamma_{c}(\clg,r^{*}\cla)$, $ ||f||_{red}=\sup\{||\Lambda_{x}(f)||: x \in \clg^{(0)}\}$. This is because the family of $*$-morphisms $\{\epsilon_{x}: x \in \mathcal{G}^{(0)}$\} separates elements of $B$. For the same reason, it follows that if $D$ is a dense subset of $\clg^{(0)}$, $||f||_{red}:=\sup\{||\Lambda_{x}(f)||: x \in D\}$ for $f \in \Gamma_{c}(\clg,r^{*}\cla)$.
\item[(2)] Let $\gamma \in \clg$ be such that $s(\gamma)=x$ and $r(\gamma)=y$. Then $||\Lambda_{x}(f)||=||\Lambda_{y}(f)||$ for $f \in \Gamma_{c}(\clg,r^{*}\cla)$.
\item[(3)] Thus if $x_{0} \in \clg^{(0)}$ is such that its $\clg$-orbit i.e. \[\{y \in \clg^{(0)}:\textrm{ there exists $\gamma \in \clg$ such that $s(\gamma)=x$ and $r(\gamma)=y$}\}\] is dense in $\clg^{(0)}$ then $\Lambda_{x_{0}}$ is a faithful representation of $\cla \rtimes_{red} \clg$. 
\end{enumerate}
 
\begin{rmrk}
\label{some remarks}
We need the following facts whose proofs we invite the reader to supply.
\begin{enumerate}
\item[(1)] If $\mathcal{A}:=\clg^{(0)} \times \mathbb{C}$ and if $\alpha_{\gamma}$ is the identity map for $\gamma \in \clg$, then $\cla \rtimes_{red} \clg$ is isomorphic to $C_{red}^{*}(\clg)$.
\item[(2)] Let $\clg:=X \rtimes G$ be a transformation groupoid where the  group $G$ acts on $X$ on the right. We assume of course that $G$ is locally compact and $X$ is locally compact and Hausdorff. Suppose $(\cla,\clg, \alpha)$ be a groupoid dynamical system. Then the group $G$ acts on the $C^{*}$-algebra  $\Gamma_{0}(X,\cla)$ where the action $\widetilde{\alpha}$ is as follows : For $g \in G$, and $f \in \Gamma_{0}(X,\cla)$, let $\widetilde{\alpha}_{g}(f)(x):=\alpha_{(x,g)}(f(x.g))$ for $x \in X$. Then $\mathcal{A} \rtimes_{red} \clg$ is isomorphic to the reduced crossed product $\Gamma_{0}(X,\cla) \rtimes_{\widetilde{\alpha},red} G$.
\end{enumerate} 
 \end{rmrk}

\section{Wiener-Hopf groupoid}

Let $G$ be a locally compact, second countable, Hausdorff topological group.  Let $P \subset G$ be a closed subsemigroup containing the identity element $e$. We say that 
\begin{enumerate}
\item[(a)] $P$ is a right Ore subsemigroup of $G$ if $PP^{-1}=G$, and
\item[(b)] $P$ is solid if $Int(P)$ is dense in $P$. 
\end{enumerate}
We say that $P$ is  a left Ore subsemigroup of $G$ if  $P^{-1}P=G$. For the rest of this article, we reserve the letters $P$ and $G$ to mean that $P$ is a closed, solid right Ore subsemigroup of a locally compact group $G$ containing the identity element  $e$ of $G$. We denote the left Haar measure on $G$ by $dg$ and the modular function by $\Delta$.

For $a,b \in P$, we say $a \leq b$ ($a<b$)  if and only if $a^{-1}b \in P$ ($a^{-1}b \in Int(P)$). If $P$ is right Ore then the preorder $\leq$ on $P$ is directed i.e. given $a,b \in P$ there exists $c \in P$ such that $a,b \leq c$.  

\begin{rmrk}
Note that $Int(P)$ is closed under right and left multiplication by $P$. Thus $G:=Int(P)Int(P)^{-1}$ if $P$ is a solid right Ore subsemigroup of $G$. Note that $(Int(P),<)$ is directed.
\end{rmrk}

Let $X$ be a second countable, locally compact Hausdorff topological space. Let $\mathcal{C}(X)$ be the set of closed subsets of $X$. We recall here the essential facts regarding the Fell topology (\cite{Fell}) on $\mathcal{C}(X)$.  For a compact subset $K \subset X$ and an open subset $U \subset X$, let \[
\mathcal{U}(K,U):=\{A  \in \mathcal{C}(X): A \cap K=\emptyset ~\textrm{and~} A \cap U \neq \emptyset\}.\]
The family $\{\mathcal{U}(K,U): K \textrm{~compact~}, U \textrm{~open}\}$ forms a sub base for a topology on $\mathcal{C}(X)$, called the Fell topology on $\mathcal{C}(X)$.  We refer the reader to Chapter 5 of  \cite{Beer}  for proofs.  
Endowed with the Fell topology, $\mathcal{C}(X)$ is compact(Thm. 5.1.3,  \cite{Beer}).  Since we have assumed that $X$ is locally compact and second countable, $\mathcal{C}(X)$ with the Fell topology is metrizable (Thm. 5.1.5, \cite{Beer}). 

Let $(A_{n})$ be a sequence of closed subsets of $X$. Define the sets  $\liminf A_{n}$ and $\limsup A_{n}$ as follows:  For $x \in X$, 
\begin{itemize}
\item  $x \in \liminf A_{n}$ if and only if there exists $N \in \mathbb{N}$ and  a sequence $(x_{n})_{n \geq N}$ in $X$ such that $x_{n} \in A_{n}$ and $x_{n} \to x$, and
\item  $x \in \limsup A_{n}$ if and only if there exists a sequence $(x_{n_k})$ in $X$, where $(n_k)$ is a subsequence, such that $x_{n_k} \in A_{n_k}$ and $x_{n_k} \to x$.  
\end{itemize}

Note that $\liminf A_{n} \subset \limsup A_{n}$ and $\liminf A_{n}$ and $\limsup A_{n}$ are both closed.  

We need the following  criterion (Thm. 5.2.6, \cite{Beer}) for the convergence of the sequence $(A_{n})$ in $\mathcal{C}(X)$.
Let $(A_{n})$ be a sequence of closed subsets of $X$. The sequence $(A_{n})$ converges in $\mathcal{C}(X)$ if and only
if $\limsup A_{n}=\liminf A_{n}$. Moreover if $\limsup A_n=\liminf A_{n}$ then $A_{n}$ converges to $A:=\limsup A_{n}=\liminf A_{n}$.

\begin{rmrk}
The  Fell topology on $\mathcal{C}(X)$, in the context of Wiener-Hopf $C^{*}$-algebras, was first exploited by Hilgert and Neeb in \cite{Hilgert_Neeb}. They introduce a topology on $\mathcal{C}(X)$ in a slightly different way by embedding $\mathcal{C}(X)$ into the space of non-empty compact subsets , endowed with the Vietoris topology, of the one point compactification of $X$. But Lemma I.4 of \cite{Hilgert_Neeb} and Thm. 5.2.6. of \cite{Beer} implies that the topology on $\mathcal{C}(X)$ used in \cite{Hilgert_Neeb} agrees with the Fell topology on $\mathcal{C}(X)$.  For similarities and differences between the Vietoris topology and the Fell topology, we refer the reader to the book \cite{Beer}.
For the definition of the Vietoris topology, we refer the reader to Chapter 2 of \cite{Beer}.
\end{rmrk}

Let $\mathcal{C}(G)$ be the space of closed subsets of $G$ endowed with the Fell topology. The group $G$ acts on $\mathcal{C}(G)$ on the right by right multiplication i.e. for $X \in \mathcal{C}(G)$ and $g \in G$, $X.g=Xg$. It is easy to verify that the action of $G$ on $\mathcal{C}(G)$ is continuous  (See Prop. I.5, \cite{Hilgert_Neeb}).

Consider the transformation groupoid $\mathcal{C}(G) \rtimes G$. Let 
\begin{equation}{\label{units}}
\Omega:= \overline{\{P^{-1}a: a \in P\}}.
\end{equation}
Then $\Omega$, being  the closure of $\{P^{-1}a: a \in P\}$ in $\mathcal{C}(G)$, is compact and is  invariant under the right action of $P$. We call $\Omega$ the order compactification of $P$. Let $\mathcal{G}:=\Omega \rtimes P$ be the semi-direct product groupoid (See \cite{Jean_Sundar} for its construction). In this case,  $\mathcal{G}$ is just the restricted groupoid $\mathcal{C}(G)\rtimes G|_{\Omega}$ i.e. $\mathcal{G}:=\{(X,g) \in \Omega \times G : X.g \in \Omega\}$ where the groupoid operations are given by the formula:
\begin{align*}
(X,g)(X.g,h)&=(X,gh) \\
(X,g)^{-1}&=(X.g,g^{-1})
\end{align*}
We call the semi-direct product groupoid $\Omega \rtimes P$ as the Wiener-Hopf groupoid associated to the pair $(P,G)$.

\begin{rmrk}
\label{facts from Ren_Sundar}
We recall some of the important results from \cite{Jean_Sundar} needed in this paper.
\begin{enumerate}
\item[(1)] The Wiener-Hopf groupoid has a left Haar system. For $X \in \Omega$, let \[Q_{X}:=\{g \in G: X.g \in \Omega\}.\] Then $\mathcal{G}^{(X)}={X} \times Q_{X}$ for $X \in \Omega$. With this identification, a left Haar system $(\lambda^{(X)})_{X\in \Omega}$ on $\clg$ is given by $\lambda^{(X)}:=1_{Q_{X}}dg$. We use only this left Haar system in this paper. Also for $X \in \Omega$, $Q_{X}=X^{-1}$ (See Proposition 5.1 of \cite{Jean_Sundar}).

\item[(2)] The crucial fact in proving (1) is that if $(X_{n})$ is a sequence in $\Omega$ converging to $X$ then $1_{X_{n}}\to 1_{X}$ a.e. (See the proof of the implication $(4)\Rightarrow (5)$ of Theorem 4.3 of \cite{Jean_Sundar}). 

\item[(3)] For  $X \in \Omega$, $Int(X)$ is dense in $X$. (See Lemma 4.1 of \cite{Jean_Sundar}) 
\end{enumerate}
\end{rmrk}
Let $\widetilde{\Omega}:=\bigcup_{a \in P} \Omega a^{-1}$. Note that $\widetilde{\Omega}$ is $G$-invariant. For  $\Omega$ is $P$-invariant and $PP^{-1}=G$. Let $\Omega_{0}:=\Omega Int(P)$. Then
\begin{equation}
\label{Omega0}
\Omega_{0}=\{X \in \Omega: X \cap Int(P) \neq \emptyset \}.
\end{equation}
(See the proof of Proposition 5.1 in \cite{Jean_Sundar}). Hence $\Omega_{0}$ is open in $\Omega$.

\begin{lmma}
\label{equivalence}
We have the following.
\begin{enumerate}
\item[(1)] The set $\Omega_{0}$ is open in $\widetilde{\Omega}$.
\item[(2)]  $\widetilde{\Omega}=\bigcup_{a \in Int(P)}\Omega_{0}a^{-1}=\bigcup_{a \in P}\Omega a^{-1}$.
\item[(3)] The space $\widetilde{\Omega}$, equipped with the subspace topology, is locally compact. 
\item[(4)] Let $A \in \widetilde{\Omega}$ and $g \in G$ be given. Then $Ag \in \Omega$ if and only $g^{-1} \in A$.
\end{enumerate}
The Wiener-Hopf groupoid $\Omega \rtimes P$ is equivalent to $\widetilde{\Omega} \rtimes G$.
\end{lmma}
\textit{Proof.} We claim that $\Omega_{0}:=\{X \in \widetilde{\Omega}: X \cap Int(P) \neq \emptyset\}$. Clearly the inclusion $\subset$ holds. Now let $X \in \widetilde{\Omega}$ be such that $X \cap Int(P) \neq \emptyset$. Write $X:=Ya^{-1}$ for $Y \in \Omega$ and $a \in P$. Choose a sequence $b_{n} \in P$ such that $P^{-1}b_{n} \to Y$. Then $P^{-1}b_{n}a^{-1} \to X$. Since  $\{Z \in \mathcal{C}(G):  Z \cap Int(P) \neq \emptyset\}$ is open, it follows that $P^{-1}b_{n}a^{-1}$ intersects $Int(P)$ eventually. But $PInt(P) \subset Int(P)$. Hence $b_{n}a^{-1} \in Int(P)$ eventually and it follows that $X$, being the limit of $P^{-1}b_{n}a^{-1}$, belongs to $\Omega$. From Equation \ref{Omega0}, it follows that $X \in \Omega_{0}$. Now the openness of $\Omega_{0}$ in $\widetilde{\Omega}$ follows from the equality $\Omega_{0}:=\{X \in \widetilde{\Omega}: X \cap Int(P) \neq \emptyset\}$. This proves $(1)$.

The proof of $(2)$ follows from the fact that $Int(P)P \subset Int(P)$ and $Int(P)$ is non-empty. Since $G$ acts continuously on $\widetilde{\Omega}$, it follows that for $a \in P$, $\Omega_{0}a^{-1}$ is open in $\widetilde{\Omega}$ and its closure in $\widetilde{\Omega}$ is the compact set $\Omega a^{-1}$. The local compactness of $\widetilde{\Omega}$ follows from this. This proves $(3)$.

Let $A \in \widetilde{\Omega}$ and $g \in G$ be given. Choose $B \in \Omega$ and $a \in P$ such that $A=Ba^{-1}$.  By (1) of Remark \ref{facts from Ren_Sundar}, we have that
$Ag \in \Omega \Leftrightarrow Ba^{-1}g \in \Omega \Leftrightarrow g^{-1}a \in B \Leftrightarrow g^{-1}\in Ba^{-1}=A.$ This proves $(4)$.

The equivalence betweem $\Omega \rtimes P$ and $\widetilde{\Omega} \rtimes G$ follows from Prop. 6.4 and Thm 6.2 of \cite{Jean_Sundar}. This completes the proof. \hfill $\Box$

 \section{The Toeplitz algebra associated to Endomorphisms of  $P$}
 
 Let $A$ be a $C^{*}$-algebra and let $End(A)$ be the set of $*$-algebra homomorphisms on $A$. By a strongly continuous left action of $P$ on $A$, we mean  a map $\alpha:P \to End(A)$ such that 
 \begin{enumerate}
\item for $a,b \in P$, $\alpha_{ab}=\alpha_{a}\alpha_{b}$,
\item $\alpha_{e}=id_{A}$ , and
\item for $x \in A$, the map $P \ni a \to \alpha_{a}(x) \in A$ is continuous. 
  \end{enumerate}
  We say the action $\alpha$ is injective/surjective if for every $a \in P$, $\alpha_{a}$ is injective/surjective. If $A$ is unital, we say $\alpha$ is unital if $\alpha_{a}$ is unital for $a \in P$. Since we will consider only strongly continuous actions of $P$, we will drop the adjective strongly continuous from now on.
  
\begin{xmpl}
Let $\clh$ be a separable Hilbert space and $\cla(\clh)$ be the CAR algebra on $\clh$ i.e. the universal unital $C^{*}$-algebra generated by symbols $\{a(\xi):\xi \in \clh\}$ such that 
\begin{align*}
a(\xi)a(\eta)+a(\eta)a(\xi)&=0 \\
a(\xi)^*a(\eta)+a(\eta)a(\xi)^{*}&=<\eta,\xi> \\
a(\lambda \xi+ \mu \eta)& = \lambda a(\xi)+ \mu a(\eta)
\end{align*}
for $\xi, \eta \in \clh$ and $\lambda,\mu \in \mathbb{C}$. 

Denote the set of isometries on $\clh$ by $\mathcal{V}(\clh)$. Suppose $W:P \to \mathcal{V}(\clh)$ is a strongly continuous semigroup homomorphism. For $a \in P$,  define $\alpha_{a}: \cla(\clh) \to \cla(\clh)$ by $\alpha_{a}(a(\xi))=a(W_{a}\xi)$. Then $\alpha: P \to End(\cla(\clh))$ is a unital injective action. The action is injective. For, the CAR algebra is simple. We can take $\clh$ to be $L^{2}(P)$ and $W:P \to \mathcal{V}(L^{2}(P))$ to be the left regular representation.
\end{xmpl}  

\begin{xmpl}
Let $X$ be a compact space and suppose $P$ acts continuously on the right injectively on $X$. Then for $a \in P$, let $\alpha_{a}:C(X) \to C(X)$ be defined by \[\alpha_{a}(f)(x):=f(xa)\]
for $f \in C(X)$. Then $\alpha:P \to End(C(X))$ is a surjective action.
\end{xmpl}
 
  We will consider the Hilbert space $L^{2}(P)$ as a closed subspace of $L^{2}(G)$ by extending an element of $L^{2}(P)$ to $L^{2}(G)$ by declaring its value outside $P$ to be zero. For $g \in G$, let $w_{g}:L^{2}(P) \to L^{2}(P)$ be defined by the formula : For $\xi \in L^{2}(P)$, let 
 \[
 w_{g}\xi(a)=\xi(ag^{-1})\Delta(g)^{-\frac{1}{2}} \textrm{ for $a \in P$}.
 \]
 For $a \in P$, let $v_{a}:=w_{a}$. Then 
 \begin{itemize}
 \item for $a \in P$, $v_{a}$ is an isometry,
 \item for $a,b \in P$, $v_{ab}=v_{b}v_{a}$,
  \item if $g=ab^{-1}$ with $a,b \in P$ then $w_{g}=v_{b}^{*}v_{a}$, and
  \item the map $G \ni g \to w_{g} \in B(L^{2}(P))$ is strongly continuous.
  \end{itemize}
 
 Let $A$ be a $C^{*}$-algebra and $\alpha:P \to End(A)$ be an action. Consider the Hilbert $A$-module $E:=L^{2}(P) \otimes A$. For $g \in G$, let $W_{g}:=w_{g} \otimes 1$ and let $V_{a}:=v_{a} \otimes 1$ for $a \in P$. The map $G \in g \to W_{g} \in \mathcal{L}_{A}(E)$ is strictly continuous. For $x \in A$, let $\pi(x) \in \mathcal{L}_{A}(E)$ be the operator defined on the dense subspace $C_{c}(P,A) \subset E$ as follows : For $\xi \in C_{c}(P,A)$, \[
 (\pi(x)\xi)(a):=\alpha_{a}(x)\xi(a).
 \] 
 The covariance condition $V_{a}^{*}\pi(x)V_{a}=\pi(\alpha_{a}(x))$ holds for every $a \in P$ and $x \in A$.
 
  To verify this, let $\rho: A \to B(\clh)$ be a faithful non-degenerate representation on a Hilbert space $\clh$. Consider the Hilbert space $\widetilde{\clh}:=L^{2}(P,\clh)=L^{2}(P) \otimes \clh$. Let $\widetilde{\rho}:A \to B(\widetilde{\clh})$ be the representation defined by the formula: For $x \in A$, $\xi \in \widetilde{\clh}$, let 
 \[
  (\widetilde{\rho}(x)\xi)(a):=\rho(\alpha_{a}(x))\xi(a) \textrm{ for $a \in P$}.
 \]
 For $g \in G$, let $\widetilde{W_{g}} \in B(\widetilde{\clh})$ be defined by 
 $\widetilde{W_{g}}\xi(b):=\xi(bg^{-1})\Delta(g)^{-\frac{1}{2}}$ for $\xi \in \widetilde{\clh}$ and $a \in P$. For $a \in P$, let $\widetilde{V_{a}}=W_{a}$. Then $\widetilde{V_{a}}^{*}=W_{a^{-1}}$ for $a \in P$. Then it is easily verifiable that \[\widetilde{V_{a}}^{*}\widetilde{\rho}(x)\widetilde{V_{a}}=\widetilde{\rho}(\alpha_{a}(x))\] for $x \in A$ and $a \in P$.
 
 Consider the internal tensor product $E \otimes_{\rho} \clh$ which is a Hilbert space. The map $\mathcal{L}_{A}(E) \ni T \to T \otimes 1 \in B(E \otimes_{\rho} \clh)$ is injective. Define $U:(L^{2}(P) \otimes A) \otimes_{\rho} \clh \to \widetilde{\clh}$ by $U(( \xi \otimes x) \otimes \eta) =\xi \otimes \rho(x)\eta $ for $\xi \in L^{2}(P)$, $x \in A$ and $\eta \in \clh$.  Then $U$ extends to a unitary from $E \otimes_{\rho} \clh$ to $\widetilde{\clh}$. We leave it to the reader to verify that for $g \in G$, $x \in A$,
 \begin{align*}
 U(\pi(x) \otimes 1)U^{*}&= \widetilde{\rho}(x) \\
 U(W_{g} \otimes 1)U^{*}&= \widetilde{W}_{g}.
 \end{align*}
 Hence it follows that $V_{a}^{*}\pi(x)V_{a}=\pi(\alpha_{a}(x))$ for $x \in A$ and $a \in P$. Henceforth, we will simply denote $\pi(x)$ by $x$ for $x \in A$. Also note that the relation $V_{a}^{*}xV_{a}=\alpha_{a}(x)$ and the fact that $V_{a}V_{a}^{*}$ commutes with $x$ implies that $V_{a}^{*}x=\alpha_{a}(x)V_{a}^{*}$ for $a \in P$ and $x \in A$.

 For $f \in C_{c}(G,A)$, let 
 \[
 W_{f}:=\int f(g)W_{g} dg \in \mathcal{L}_{A}(E).
 \]
  Since the map $g \to W_{g}$ is strictly continuous, it follows that the above integral exists in $\mathcal{L}_{A}(E)$. For $f \in C_{c}(G,A)$, let us call $W_{f}$ the Wiener-Hopf operator with symbol $f$.
  
   \textbf{\emph{The Toeplitz algebra}}  denoted $\mathfrak{T}(A,P,\alpha)$ is defined as the $C^{*}$-subalgebra of $\mathcal{L}_{A}(E)$ generated by $\{W_{f}: f \in C_{c}(G,A)\}$. We will denote the $C^{*}$-subalgebra of $\mathcal{L}_{A}(E)$ generated by $\{W_{f}: f \in C_{c}(Int(P),A)\}$ by $A \rtimes_{red,\alpha} P$ and call it the \textbf{\emph{reduced crossed product}} $A$ by $P$. We will denote $A \rtimes_{red,\alpha} P$ by $A \rtimes_{red} P$ if the action $\alpha$ is clear.
 
 \begin{rmrk}
 Alternatively, $\mathfrak{T}(A,P,\alpha)$ can defined (by choosing a faithful non-degenerate representation $\rho$ of $A$ on $\clh$) as the $C^{*}$-subalgebra of $B(\widetilde{\clh})$ generated by the integrals \[\Big \{\int \widetilde{\rho}(f(g)) \widetilde{W}_{g} dg: f \in C_{c}(G,A)\Big \}\] where $\widetilde{\clh}$, $\widetilde{\rho}$ and $\widetilde{W}_{g}$ are defined as in the preceeding paragraphs.
 \end{rmrk}
 
  The Toeplitz algebra $\mathfrak{T}(A,P,\alpha)$ was studied by Khoskham and Skandalis in \cite{KS97} for the additive semigroups $\mathbb{N}$ and $[0,\infty)$. In this paper, we describe $\mathfrak{T}(A,P,\alpha)$ as a groupoid crossed product when either  the action $\alpha$ is unital and injective or when the action $\alpha$ is surjective. Moreover we show in both cases that the two algebras $\mathfrak{T}(A,P,\alpha)$ and $A \rtimes_{red} P$ coincide. When $A=\mathbb{C}$, then $\mathfrak{T}(\mathbb{C},P, \alpha)$ is the usual Wiener-Hopf algebra whose study, from the groupoid point of view, was initiated in \cite{Renault_Muhly}.

 We end this section with a lemma which is essential in proving $\mathfrak{T}(A,P,\alpha)=A \rtimes_{red} P$.
  \textit{Notation:} For $f \in C_{c}(G)$ and $a \in G$, let $L_{a}(f) \in C_{c}(G)$ be defined by the formula $L_{a}(f)(x)=f(a^{-1}x)$ for $x \in G$.
 
 \begin{lmma}
 \label{convolution_density}
 Let \[\mathcal{F}:= \{\phi*\psi : \phi \in C_{c}(Int(P)) \textrm{~and~} \psi \in C_{c}(Int(P^{-1}))\}.\] Then the linear span of $\mathcal{F}$ is dense in $C_{c}(G)$ in the inductive limit topology. 
 \end{lmma}
 \textit{Proof.} The proof is an imitation of the proof of Lemma 2.3 of \cite{Jean_Sundar} where the density of $span(\mathcal{F})$ is proved in  the $L^{1}$-norm. Let $f \in C_{c}(G)$ be given and $K:=supp(f)$. Since  $G=Int(P)Int(P)^{-1}$, it follows that there exists $a_{1},a_{2},\cdots,a_{m} \in Int(P)$ such that $K \subset \displaystyle \bigcup_{i=1}^{m} a_{i}Int(P)^{-1}$. By choosing a partition of unity, we can write $f=\sum_{i=1}^{n}f_{i}$ with $f_{i} \in C_{c}(G)$ and $supp(f_{i}) \subset a_{i}Int(P)^{-1}$. For $i=1,2,\cdots, m$, let $\widetilde{f_{i}}:=L_{a_{i}^{-1}}f_{i}$. Note that $supp(\widetilde{f_{i}}) \subset Int(P)^{-1}$ for every $i$.
 
 Note that the intersection $U:=\displaystyle \bigcap_{i=1}^{m}Int(P)^{-1}a_{i}$ contains the identity element $e$ and hence $U \cap P$ is non-empty. Since $Int(P)$ is dense in $P$, it follows that $U \cap Int(P) \neq \emptyset$. Choose $b \in U \cap Int(P)$. Then $b^{-1}Int(P)$ is an open set containing $e$. 
 
 Now choose a countable base $\{U_{n}: n \in \mathbb{N}\}$, which we can assume to be decreasing, at $e$ such that $U_{n} \subset b^{-1}Int(P)$. We can also assume that each $U_{n}$ has compact closure.  Choose $\phi_{n} \in C_{c}(G)$ such that $\phi_{n} \geq 0$, $supp(\phi_{n}) \subset U_{n}$ and $\int \phi_{n}(t)dt =1$. Then $\{\phi_{n}\}_{n \in \mathbb{N}}$ forms an approximate identity for the convolution product when $C_{c}(G)$ is given the inductive limit topology.  
  Note that $L_{a_{i}}\phi_{n}$ is supported in $a_{i}b^{-1}Int(P)$ for every $i$. 
 But the choice of $b$ implies that $a_{i}b^{-1} \in Int(P)$. Hence for every $i$, $supp(L_{a_{i}}\phi_{n}) \subset Int(P)$. 
 
 Note that $\displaystyle \mathcal{F} \ni \sum_{i=1}^{m}L_{a_{i}}\phi_{n}*\widetilde{f_{i}} = \sum_{i=1}^{m}L_{a_{i}}(\phi_{n}*\widetilde{f_{i}})$. Since left translations are continuous linear operators in the inductive limit topology, it follows that as $n$ tends to infinity, $\sum_{i=1}^{m}L_{a_{i}}\phi_{n}*\widetilde{f_{i}}$ converges, in the inductive limit topology, to $\sum_{i=1}^{m}L_{a_{i}}\widetilde{f_{i}}=f$.  This completes the proof. \hfill $\Box$

 \section{The injective case}
 
 Throughout this section, we assume that $A$ is a unital separable  $C^{*}$-algebra and the strongly continuous action $\alpha:P \to End(A)$  is unital and injective. We also assume that $P$ is left Ore i.e. $P^{-1}P=G$. Note that in this case $G=Int(P)^{-1}Int(P)$.
 
 First we dilate the action of $P$ on $A$ to an action of $G$ by adapting Laca's dilation theorem (Thm 2.1.1, \cite{Laca2000}). Note that Laca's dilation theorem is for discrete Ore semigroups and we adapt his argument to the locally compact case.
 
 Consider for the moment, by forgetting the topology on $G$ and $P$, that $G$ is a discrete group and $P$ is a discrete left Ore subsemigroup of $G$.  By Laca's dilation theorem ( Thm 2.1.1, \cite{Laca2000}), there exists a unital $C^{*}$-algebra $B$, a left $G$-action $\widetilde{\alpha}:G \to Aut(B)$ and an  unital embedding $j: A \to B$ such that 
 \begin{enumerate}
 \item[(1)] for $a \in P$ and $x \in A$, $j(\alpha_{a}(x))=\widetilde{\alpha}_{a}(j(x))$, and
 \item[(2)] the union $\bigcup_{a \in P}\widetilde{\alpha}_{a}^{-1}(j(A))$ is dense in $B$.
 \end{enumerate}
Moreover the $G$-$C^{*}$-algebra $B$, together with the embedding $j:A \to B$, is unique up to a $G$-equivariant isomorphism. We will fix such a $G$-$C^{*}$-algebra $B$ and will identify $A$ with a unital $C^{*}$-subalgebra of $B$. For $g \in G$, we will simply denote $\widetilde{\alpha}_{g}$ by $\alpha_{g}$. 
 
\textit{Claim:} The action of $G$ on $B$ is strongly continuous. 

Since $B$ is the closure of $ \bigcup_{a \in P}\alpha_{a}^{-1}(A)$, it is enough to show that for $x \in A$ and for $a \in P$, the map $G \ni g \to \alpha_{g}\alpha_{a}^{-1}(x)=\alpha_{ga^{-1}}(x) \in B$ is continuous. The map $G \ni g \to ga^{-1} \in G$ is a homeomorphism for $a \in P$. Thus it is enough to show that for $x \in A$, the map $G \ni g \to \alpha_{g}(x) \in B$ is continuous.

Let $x \in A$ and let $(g_{n})$ be a sequence in $G$ such that $(g_{n}) \to g$. Write $g=a^{-1}b$ with $a,b \in Int(P)$. Then $g \in a^{-1}Int(P)$. Thus eventually $g_{n} \in a^{-1}Int(P)$. For large $n$, write $g_{n}=a^{-1}b_{n}$ with $b_{n} \in Int(P)$. Then $b_{n} \to b$. Since the action $\alpha:P \to End(A)$ is strongly continuous, it follows that $\alpha_{b_{n}}(x) \to \alpha_{b}(x)$. Hence  $\alpha_{g_{n}}(x)=\alpha_{a}^{-1}\alpha_{b_{n}}(x)$ converges to  $ \alpha_{a}^{-1}\alpha_{b}(x)=\alpha_{g}(x)$. This proves the claim.

Note that since $A$ is separable, the equality $B=\overline{\bigcup_{a \in P}\alpha_{a}^{-1}(A)}$ and the fact that the action is strongly continuous implies that $B$ is separable. Recall that we have assumed that $P$ is second countable. 

 Thus for the rest of the section,  \textit{we assume that $B$ is a unital separable $C^{*}$-algebra with a strongly continuous action $\alpha: G \to Aut(B)$ and $A \subset B$ is a  unital (i.e. $1_{B} \in A$) $C^{*}$-subalgebra  which is  $P$-invariant. }

 Let $\mathcal{G}:=\Omega \rtimes P$ be the Wiener-Hopf groupoid associated to $P$. Recall that $\Omega$ is the closure of $\{P^{-1}a: a \in P\}$ in the space of closed sets of $G$ equipped with the Fell topology.  For $X \in \Omega$, let $A_{X}:=C^{*}\{\alpha_{g}^{-1}(x): g \in X, x \in A \}$ i.e. $A_{X}$ is the $C^{*}$-subalgebra of $B$ generated by $\bigcup_{g \in X}\alpha_{g}^{-1}(A)$.
 
  Observe the following.
 \begin{itemize}
 \item For $a \in P$, $A_{P^{-1}a}=\alpha_{a}^{-1}(A)$.
 \item Let $g \in G$ and $X,Y \in \Omega$ be  such that $Y=Xg$. Then the map $\alpha_{g}$ maps $A_{Y} $ onto $A_{X}$ and is an isomorphism from $A_{Y}$ onto $A_{X}$.
 \item For $X \in \Omega$, $A \subset A_{X}$. This is because $e \in X$ for $X \in \Omega$.
 \end{itemize} 
 
 Let  \begin{equation}
 \label{bundle-injective}
  \cla:=\{(X,x): X \in \Omega, x \in A_{X}\}
  \end{equation} 
  We equip $\cla$ with the subspace topology induced from the product topology on $ \Omega \times B$. Denote the projection map $\cla \in (X,x) \to X \in  \Omega$ by $p$. Note that  $p^{-1}(X)=A_{X}$ for $X \in \Omega$. For $(X,g) \in \Omega \rtimes P$, let $\alpha_{(X,g)}:A_{X.g} \to A_{X}$ be the map $\alpha_{g}$ restricted to $A_{X.g}$. Thus the Wiener-Hopf groupoid $\clg$ acts on the bundle $\cla$ via the maps $\{\alpha_{(X,g)}:(X,g) \in \clg\}$. We  denote both the action of $\clg$ on $\cla$ and the action of the group $G$ on $B$ by the same letter $\alpha$. 

We will identify a section of the bundle $\cla$ with a map $F:\Omega \to B$ such that $F(X) \in A_X$ for every $X \in \Omega$.   
 
 \begin{lmma}
 \label{density}
 Let $x \in A$ and $f \in C_{c}(G)$ be given. Define $F_{x,f}:\Omega \to \cla $ by the formula
 \[
 F_{x,f}(X):=\int \alpha_{g}^{-1}(x)f(g)1_{X}(g)dg.
 \]
 Then $F_{x,f}$ is a continuous section (i.e. $p \circ F_{x,f}=id$). 
 
 Also for every $X \in \Omega$, the $*$-subalgebra generated by $\{F_{x,f}(X):x \in A, f \in C_{c}(G)\}$ is dense in $A_{X}$.
 \end{lmma}
 \textit{Proof.} Let $x \in A$ and $f \in C_{c}(G)$ be given. By definition, it follows that $F_{x,f}(X) \in A_{X}$ for $X \in \Omega$. Let $(X_{n}) \to X $ in $\Omega$. By Remark \ref{facts from Ren_Sundar}, we have $1_{X_{n}} \to 1_{X}$ a.e. Then by the dominated convergence theorem, it follows that $F_{x,f}(X_{n}) \to F_{x,f}(X)$. Hence $F_{x,f}$ is continuous.
 
Let $X \in \Omega$ be given. Since $Int(X)$ is dense in $X$ (Remark \ref{facts from Ren_Sundar}), it follows that  $\displaystyle A_{X}$ is the $C^{*}$-subalgebra generated by $\{\alpha_{g}^{-1}(x): g \in Int(X),x \in A\}$. Now let $g_{0} \in Int(X)$, $x \in A$ and $\epsilon >0$ be given.  Since the action $\alpha:G \to Aut(B)$ is strongly continuous, there exists an open set $U \subset Int(X)$ containing $g_{0}$ such that $||\alpha_{g}^{-1}(x)-\alpha_{g_{0}}^{-1}(x)|| \leq \epsilon$ for $g \in U$. Choose $f \in C_{c}(G)$ such that $ f \geq 0$, $supp(f) \subset U$ and $\int f(g) dg=1$. Now observe that 
  \begin{align*}
 ||F_{x,f}(X)-\alpha_{g_{0}}^{-1}(x) || & = \Big| \Big| \int_{g \in U} \alpha_{g}^{-1}(x)f(g)1_{X}(g)dg - \int_{g \in U}\alpha_{g_{0}}^{-1}(x)f(g)dg ~\Big |\Big| \\
 & =\Big|\Big| \int_{g \in U} (\alpha_{g}^{-1}(x)-\alpha_{g_{0}}^{-1}(x))f(g)dg ~\Big|\Big| \\
 & \leq \int_{g \in U} || \alpha_{g}^{-1}(x)-\alpha_{g_{0}}^{-1}(x)|| f(g) dg \\
 & \leq \epsilon 
  \end{align*} 
 Thus the $*$-subalgebra generated by $\{F_{x,f}(X): x \in A, f \in C_{c}(G)\}$ is dense in $A_{X}$. This completes the proof. \hfill $\Box$
 
 Let $\mathcal{E}$ be the subalgebra generated by $\{F_{x,f}: x \in A, f \in C_{c}(G)\}$ in $C(\Omega,B)$ where $C(\Omega,B)$ denotes the algebra of continuous functions from $\Omega$ to $B$. Note that $\mathcal{E}$ is $*$-closed.
 
 \begin{ppsn}
 \label{bundle}
We have the following.
\begin{enumerate}
\item[(1)] The map $p: \cla \to \Omega$ is an open surjection.
\item[(2)] The triple $(\cla,p,\Omega)$ is an upper semi-continuous bundle of $C^{*}$-algebras.
\item[(3)] The linear span of $\{\psi F: F \in \mathcal{E}, \psi \in C(\Omega)\}$ is dense in $\Gamma(\Omega,\cla)$.
\item[(4)] The triple $(\cla,\clg,\alpha)$ is a groupoid dynamical system. 
\end{enumerate}
 \end{ppsn}
 \textit{Proof.} Let $U \subset \cla$ be an open set and let $V:=p(U)$. Consider a point $X \in V$. Then there exists $a \in A_{X}$ such that $(X,a) \in U$. By Lemma \ref{density}, there exists $F \in \mathcal{E}$ such that $F(X) \in U$. Since $F$ is continuous, there exists an open set $W \subset \Omega$ containing $X$ such that $F(Y) \in U$ for $Y \in W$. But $p \circ F=id$. Thus $W \subset p(U)$ and contains $X$. This shows that $p(U)$ is open. This proves $(1)$.

The other conditions needed to show $(2)$ follows trivially from the fact that $\Omega \times B$ is an upper semi-continuous (in fact continuous) bundle over $\Omega$. Lemma \ref{density} and Proposition C.24 in Appendix C of \cite{Dana-Williams} implies $(3)$. The fact that $\alpha:G \to Aut(B)$ is a strongly continuous action implies $(4)$. This completes the proof. \hfill $\Box$
   
We will identify a section of $r^{*}\cla$ over $\clg$ with a map $f:\clg \to B$ such that $f(\gamma) \in A_{r(\gamma)}$ for $\gamma \in \clg$. Thus
\[\Gamma_{c}(\clg, r^{*}\cla):= \Big \{f: \clg \to B: \textrm{$f$ is continuous, compactly supported and } f(\gamma) \in A_{r(\gamma)} \Big \}.\]

  The convolution product and involution are given by the following formulas: For $\phi,\psi \in \Gamma_{c}(\clg, r^{*}\cla)$ and $(X,s) \in \clg$,
\begin{align*}
(\phi*\psi)(X,s)&=\int \phi(X,t)\alpha_{t}(\psi(Xt,t^{-1}s))1_{X}(t^{-1})dt, \textrm{ and} \\
 \phi^{*}(X,s)&=\alpha_{s}(\phi(Xs,s^{-1}))^{*}.
\end{align*}

 For $ f \in C_{c}(G,A)$, let $\widetilde{f}: \clg \to \cla$ be defined by \[\widetilde{f}(X,s):=f(s) \textrm{ for }(X,s)\in \clg.\] Since $A \subset A_{X}$ for $X \in \Omega$, it follows that $\widetilde{f} \in \Gamma_{c}(\clg, r^{*}\cla)$.  Since we have assumed $A$ is unital, we will identify $C_{c}(G)$ with a subset of $C_{c}(G,A)$.

  Let us recall Proposition 3.5 of \cite{Renault_Muhly}.  The reader is  refered to Lemma 5.4 and Theorem 5.5 of \cite{Jean_Sundar} for a few remarks concerning Proposition 3.5 of \cite{Renault_Muhly}.
  
  \begin{thm}[Prop 3.5 \cite{Renault_Muhly}]
  \label{RW_theorem}
  For $f \in C_{c}(G)$, let $\widetilde{f} \in C_{c}(\clg)$ be defined by \[\widetilde{f}(X,s)=f(s)~\textrm{for $(X,s) \in \clg$}.\] Then the $*$-algebra generated by $\{\widetilde{f}: f \in C_{c}(G)\}$ is dense in $C_{c}(\clg)$ in the topology induced by the $||~||_{I}$-norm.
    \end{thm}
 We prove an analog  of the above theorem for which a bit of preparation is needed.  
    For $a \in P$ and $\psi \in \Gamma_{c}(\clg,r^{*}\cla)$, let $R_{a}(\psi) \in \Gamma_{c}(\clg,r^{*}\cla)$ be defined by \[R_{a}(\psi)(X,s)=\alpha_{a}(\psi(X.a,a^{-1}s)).\]
  
 \begin{lmma}
 \label{continuity of R_{a}}
 Let $\psi \in \Gamma_{c}(\clg,r^{*}\cla)$ and let $(a_{n})$ be a sequence in $P$ such that $a_{n} \to a \in P$.
 Then $R_{a_{n}}(\psi) \to R_{a}(\psi)$ in the inductive limit topology.
  \end{lmma} 
 \textit{Proof.}   Let $\Omega \times K $ be a compact set containing the support of $\psi$ and $L$ be a compact subset of $P$ containing $\{a_{n}:n \in \mathbb{N}\}$. Then for every $n$, $R_{a_{n}}(\psi)$ and $R_{a}(\psi)$ are supported inside $\Omega \times LK$. If $R_{a_{n}}(\psi) \nrightarrow R_{a}(\psi)$, it follows that there exists $\epsilon >0$ and a subsequence $(a_{n_{k}},X_{n_{k}},s_{n_{k}}) \in P \times \clg$ such that 
  \begin{equation}
  \label{contradiction}
  ||\alpha_{a_{n_{k}}}\big(\psi(X_{n_{k}}a_{n_{k}},a_{n_{k}}^{-1}s_{n_{k}})\big)-\alpha_{a}\big(\psi(X_{n_{k}}a,a^{-1}s_{n_{k}})\big)|| \geq \epsilon.
  \end{equation}
  This implies in particular that $s_{n_{k}} \in LK$. Hence $(X_{n_{k}},a_{n_{k}},s_{n_{k}}) \subset \Omega \times L \times LK$ which is a compact subset of $\Omega \times P \times G$. By passing to a subsequence, if necessary, we can assume that $(X_{n_{k}},a_{n_{k}},s_{n_{k}})$ converges and let its limit be $(X,a,s)$. Since $\alpha$ is strongly continuous and $\psi$ is continuous, it follows that the sequences $\alpha_{a_{n_{k}}}\big(\psi(X_{n_{k}}a_{n_{k}},a_{n_{k}}^{-1}s_{n_{k}})\big)$ and $\alpha_{a}\big(\psi(X_{n_{k}}a,a^{-1}s_{n_{k}})\big)$ converge to $\alpha_{a}(\psi(Xa,a^{-1}s))$. Hence their difference converge to zero which is a contradiction to  Eq. \ref{contradiction}. This completes the proof. \hfill $\Box$
 
The  following corollary is an immediate consequence of the preceeding lemma. We leave the proof  to the reader.
  \begin{crlre}
 \label{continuity}
 Let $\psi \in \Gamma_{c}(\clg,r^{*}\cla)$ and $a \in P$ be given. Given $\epsilon >0$, there exists an open set $U \subset G$ containing $a$ such that  
 \[
 ||R_{b}(\psi)(X,s)-R_{a}(\psi)(X,s)|| \leq \epsilon
 \]
 for $b \in U \cap P$ and $(X,s) \in \clg$.
 \end{crlre}
  
Now we prove the analog of Theorem \ref{RW_theorem}. We use the following notation in the proof. For $\phi \in C(\Omega)$ and $f \in C_{c}(G)$, let $\phi \otimes f \in C_{c}(\clg)$ be the function defined by  $(\phi \otimes f)(X,s)=\phi(X)f(s)$  for $(X,s) \in \clg$.  For $\phi \in C(\Omega)$ and $\psi \in \Gamma_{c}(\clg,r^{*}\cla)$,  let $\phi\cdot \psi \in \Gamma_{c}(\clg,r^{*}\cla)$ be the function defined by $(\phi\cdot\psi)(X,s)=\phi(X)\psi(X,s)$. We will identify $C_{c}(\clg)$ with a $*$-subalgebra of $\Gamma_{c}(\clg,r^{*}(\cla))$. This is because we have assumed $A$ is unital and $\alpha$ is unital. 
 
 Before we prove the main proposition, we need the following lemma which is the analog of Lemma 5.4 in \cite{Jean_Sundar}.  In \cite{Jean_Sundar}, the proof is provided for the case when $A=\mathbb{C}$ and w.r.t the inductive limit topology. Here we use the $I$-norm. Since the proof is almost identical, we omit the proof. For $\phi \in \Gamma(\Omega,\cla)$ and $f \in C_{c}(G)$, we define $\phi \otimes f \in \Gamma_{c}(\clg,r^{*}\cla)$ by the equation 
 \[
 (\phi \otimes f)(X,s)=\phi(X)f(s)
 \]
 for $(X,s) \in \clg$.
 
 \begin{lmma}
 \label{les generateurs}
 Let $\mathcal{E}^{'} \subset \Gamma(\Omega, \cla)$ be a subset which is closed under conjugation. Suppose that $\Gamma(\Omega,\cla)$ is generated by $\mathcal{E}^{'}$. Then the $*$-subalgebra generated by $\{\phi \otimes f: \phi \in \mathcal{E}^{'}, f \in C_{c}(G)\}$ is dense in $\Gamma_{c}(\clg, r^{*}\cla)$ w.r.t. the $I$-norm.
  \end{lmma}
 
 Note that if we set $\mathcal{E}_{0}:=\{\phi F_{x,f}: \phi \in C(\Omega), x \in A, f \in C_{c}(G)\}$ then $\mathcal{E}_{0}$ satisfies the hypothesis of Lemma \ref{les generateurs}.
 
 \begin{ppsn}
 \label{density in inductive limit topology}
 The $*$-algebra generated by $\{\widetilde{f}: f \in C_{c}(G,A)\}$ is dense in $\Gamma_{c}(\clg,r^{*}\cla)$ in the topology induced by the norm $||~||_{I}$.
  \end{ppsn}
  \textit{Proof.} Let $\mathcal{D}$ be the $*$-subalgebra of $\Gamma_{c}(\clg,r^{*}\cla)$ generated by $\{\widetilde{f}:f \in C_{c}(G,A)\}$ and denote the closure of $\mathcal{D}$ with respect to  $||~||_{I}$ by $\overline{\mathcal{D}}$. Let $x \in A$, $f,g \in C_{c}(G)$ be given. Define $\psi:\mathcal{G} \to B$ by the formula: For $(X,s) \in \clg$, let 
  \begin{equation}
  \label{dense collection}
  \psi(X,s):= \Big(\int \alpha_{t}(x)f(t)1_{X}(t^{-1})dt\Big)g(s).
  \end{equation}
  First we claim that that there exists a sequence $\psi_{n} \in \mathcal{D}$ such that $\psi_{n} \to \psi$ in the inductive limit topology. It is enough to prove when $x$ is positive. Thus let $x$ be positive and write $x=y^{*}y$ for some $y \in A$. 
  
 Denote the support of $g$ by $K$. Let $\epsilon >0$ be given. Choose a compact open neighbourhood $V \ni e$ such that if $s,t \in G$ and $st^{-1} \in V$ then $|g(s)-g(t)| \leq \epsilon$. Since $f$ is compactly supported, there exists $t_{1},t_{2},\cdots, t_{n} \in G$ such that $supp(f) \subset \bigcup_{i=1}^{n}Vt_{i}$. Choose a partition of unity $\{\phi_{i}\}_{i=1}^{n}$ such that $\phi_{i} \geq 0$, $supp(\phi_{i}) \subset Vt_{i}$ for every $i$ and $\sum_{i=1}^{n}\phi_{i}=1$ on $supp(f)$. For every $i$, let $f_{i},g_{i} \in C_{c}(G,A)$ be defined by the equation $f_{i}(t)=y\overline{f(t^{-1})}\phi_{i}(t^{-1})$ and $g_{i}(s)=yg(t_{i}s)$. Let $\psi_{n}:=\sum_{i=1}^{n}\widetilde{f_{i}}^{*}*\widetilde{g_{i}}$. 
  
  Note that for every $i$, $\widetilde{f_{i}}^{*}*\widetilde{g_{i}}$ is supported in $\Omega \times VK$. Hence $\psi_{n}$ is supported in $\Omega \times VK$ for every $n$. Also note that for every $i$, \[|\phi_{i}(t)(g(s)-g(t_{i}t^{-1}s))| \leq \epsilon\phi_{i}(t)\] for $t,s \in G$.  If $t \in supp(\phi_{i}) \subset Vt_{i}$ then $s(t_{i}t^{-1}s)^{-1} \in V$. Then by the choice of $V$, it follows that $|\phi_{i}(t)(g(s)-g(t_{i}t^{-1}s))| \leq \epsilon \phi_{i}(t)$ for  $t \in supp(\phi_{i})$ and $s \in G$. The inequality clearly holds when $t \notin supp(\phi_{i})$
  
  Now calculate as follows to find that for $(X,s) \in \clg$, 
  \begin{align*}
  &||~ \psi(X,s) - \psi_{n}(X,s)~||\\
  & = \Big|\Big| \sum_{i=1}^{n} \Big(\int \alpha_{t}(x)f(t)\phi_{i}(t)1_{X}(t^{-1})g(s)dt - \int \alpha_{t}(y^{*})f(t)\phi_{i}(t)\alpha_{t}(y)g(t_{i}t^{ -1}s)1_{X}(t^{-1})dt \Big) \Big|\Big| \\
  & \leq \sum_{i=1}^{n} \int \big|\big|  \alpha_{t}(x)f(t)\phi_{i}(t)\big(g(s)-g(t_{i}t^{-1}s)\big)1_{X}(t^{-1})\big|\big|dt \\
  & \leq \sum_{i=1}^{n} \int ||x|| |f(t)| \epsilon \phi_{i}(t) dt \\
    & \leq \epsilon ||x|| \int |f(t)|dt.
    \end{align*}
This proves that $\mathcal{D} \ni \psi_{n} \to \psi$ in the inductive limit topology. As a consequence, it follows that $\psi \in \overline{\mathcal{D}}$.

Let $\phi \in C(\Omega)$ be given. We claim that $\phi.\psi \in \overline{\mathcal{D}}$.  We have just shown that $\psi \in \overline{\mathcal{D}}$. Consider an element $a \in Int(P)$. Choose a countable base of compact neighbourhoods  $\{U_{n}:n \in \mathbb{N}\}$ at $a$. We can assume that $U_{n}\subset Int(P)$ and is decreasing. For every $n$, let $f_{n} \in C_{c}(G)$ be such that $f_{n} \geq 0$, $supp(f_{n}) \subset U_{n}$ and  $\int f_{n}(g)dg=1$. Now observe that 
\begin{align*}
((\phi \otimes f_{n})*\psi)(X,s)&=\phi(X)\int f_{n}(t)\alpha_{t}(\psi(X.t,t^{-1}s))1_{X}(t^{-1})dt \\
                                & = \phi(X)\int _{t \in U_{n}}f_{n}(t)\alpha_{t}(\psi(X.t,t^{-1}s))dt \\
                                &=\phi(X) \int_{t \in U_{n}}f_{n}(t)R_{t}(\psi(X,s))dt .
\end{align*}
Let $\epsilon >0$ be given. Then by Corollary \ref{continuity}, it follows that there exists an open set $U \subset Int(P)$ containing $a$ such that $||R_{t}(\psi)(X,s)-R_{a}(\psi)(X,s)|| \leq \epsilon$ for $t \in U$ and $(X,s) \in \clg$. Choose $N$ such that $U_{N} \subset U$. Note that the support of $(\phi \otimes f_{n})*\psi$ is contained in $\Omega \times \overline{U_{1}}{K}$. Now observe that for $n \geq N$, 
\begin{align*}
&||((\phi \otimes f_{n})*\psi)(X,s)-\phi(X)R_{a}(\psi)(X,s)|| \\
&=\Big|\Big|\phi(X) \int_{t \in U_{n}}f_{n}(t)(R_{t}(\psi(X,s))dt-R_{a}(\psi)(X,s))dt \Big|\Big| \\
& \leq ||\phi||_{\infty} ~\epsilon \int f_{n}(t)dt \\
 &\leq \epsilon ~||\phi||_{\infty}.
\end{align*}
Hence $(\phi \otimes f_{n})*\psi \to \phi.R_{a}(\psi)$ in the inductive limit topology. By Theorem \ref{RW_theorem}, it follows that $(\phi \otimes f_{n})*\psi \in \overline{\mathcal{D}}$ for every $n \in \mathbb{N}$. Hence for every $a \in Int(P)$, $\phi \cdot R_{a}(\psi) \in \overline{\mathcal{D}}$.

 Let $(a_n)$ be a sequence in $Int(P)$ converging to $e$ where $e$ is the identity element of $G$. Then by Lemma \ref{continuity of R_{a}}, it follows that $\phi \cdot R_{a_{n}}(\psi) \to \phi \cdot R_{e}(\psi)$ in the inductive limit topology and consequently in the $I$-norm. Thus it follows that $\phi \cdot \psi \in \overline{\mathcal{D}}$. 
 
 Let $\mathcal{E}_{0}:=\{\phi F_{x,f}: x \in A, f \in C_{c}(G)\}$. We have shown that $\mathcal{E}_{0} \otimes_{alg} C_{c}(G)$ is contained in $ \overline{\mathcal{D}}$.  Now the Proposition follows  by applying Lemma \ref{les generateurs} to $\mathcal{E}_{0}$. This completes the proof. \hfill $\Box$

The rest of the discussion is similar to the case when $A=\mathbb{C}$ ( \cite{Renault_Muhly}, \cite{Jean_Sundar}). Let $X_{0}:=P^{-1} \in \Omega$. Then the $\clg$-orbit of $X_{0}$ in $\Omega$ is dense in $\Omega$. Thus the representation $\Lambda_{X_{0}}$ is a faithful representation of $\cla \rtimes_{red} \clg$.
Also note that $\mathcal{G}^{(X_{0})}=\{X_{0}\} \times P$ and $\cla_{X_{0}}=A$. Thus $L^{2}(\clg^{(X_{0})},A_{X_{0}}) = L^{2}(P,A)$ up to an obvious identification. For $f \in C_{c}(G,A)$, let $\widehat{f} \in C_{c}(G,A)$ be defined by $\widehat{f}(g)=f(g^{-1})\Delta(g)^{-\frac{1}{2}}$. Note that $C_{c}(G,A) \ni f \to \widehat{f} \in C_{c}(G,A)$ is a linear bijection.
 
 Let $\xi \in C_{c}(P,A)$ be given. Now observe that for $b \in P$,
\begin{align*}
(\Lambda_{X_{0}}(\widetilde{f})\xi)(b)&= \int \alpha_{b}(f(b^{-1}a))\xi(a)1_{P}(a) da\\
                                          & = \int \alpha_{b}(f(u))\xi(bu)1_{P}(bu) du \\
                                          &=\int \alpha_{b}(f(u^{-1}))\Delta(u)^{-\frac{1}{2}}\xi(bu^{-1})1_{P}(bu^{-1})\Delta(u)^{-\frac{1}{2}} du \\
                                         &= \big(\int \widehat{f}(u)W_{u}\xi\big)(b) \\
                \end{align*}
Thus $\Lambda_{X_{0}}(\widetilde{f})=W_{\widehat{f}}$. 

\begin{thm}
\label{main_theorem}
The Toeplitz algebra $\mathfrak{T}(A,P,\alpha)$ is isomorphic to the reduced crossed product $\cla \rtimes_{red} \clg$.
\end{thm}               
\textit{Proof.} We have just observed that  $\Lambda_{X_{0}}$ is a faithful representation of $\cla \rtimes_{red} \clg$ on the Hilbert $A$-module $L^{2}(P,A)$ and  $\Lambda_{X_{0}}(\widetilde{f}) \in \mathfrak{T}(A,P,\alpha)$ for $f \in C_{c}(G,A)$. Hence the range of $\Lambda_{X_{0}}$ contains $\mathfrak{T}(A,P,\alpha)$. Now Proposition \ref{density in inductive limit topology} implies that the range of $\Lambda_{X_{0}}$ is contained in $\mathfrak{T}(A,P,\alpha)$. Thus the range of $\Lambda_{X_{0}}$ is $\mathfrak{T}(A,P,\alpha)$. This completes the proof. \hfill $\Box$

We end this section by proving $\mathfrak{T}(A,P,\alpha)=A \rtimes_{red} P$ (under the hypotheses of this section). 

\begin{lmma}
The $C^{*}$-algebra $A \rtimes_{red} P$ contains $\{\int f(a)V_{a}^{*}da : f \in C_{c}(Int(P),A)\}$.
\end{lmma}
\textit{Proof.}  Observe that it suffices to   prove that for $x \in A$ and $f \in C_{c}(Int(P))$, the integral $\int xf(a)V_{a}^{*}da \in A \rtimes_{red} P$. Note that for $y \in A$ and $f \in C_{c}(Int(P))$, $\int \alpha_{a}(y)f(a)V_{a}^{*} \in A \rtimes_{red} P$. For the adjoint of $\int y^{*}f(a)V_{a} da$ is $\int \alpha_{a}(y)f(a)V_{a}^{*}$.
This is due to  the equality $V_{a}^{*}y=\alpha_{a}(y)V_{a}^{*}$ for $a \in P$.

Let $x \in A$ and $f \in C_{c}(Int(P))$ be given. We claim that for every $c \in Int(P)$, $\int \alpha_{c}(x)L_{c}(f)(a)V_{a}^{*}da \in A \rtimes_{red} P$.

Let $(U_{n})$ be a sequence of open sets containing the identity element $e$ such that $\{U_{n}: n \in \mathbb{N}\}$ is a neighbourhood base at $e$ and $\bigcap_{n}U_{n}=\{e\}$. We can also asumme that $\{U_{n}\}$ is decreasing. Let $c \in Int(P)$ be given. Choose $\phi_{n} \in C_{c}(G)$ such that $\phi_{n} \geq 0$, $\int \phi_{n}=1$ and $supp(\phi_{n}) \subset U_{n}$.
Let $c \in Int(P)$ be given. Define for $n \in \mathbb{N}$, $\psi_{n}:=L_{c}(\phi_{n})$. Then $supp(\psi_{n}) \subset cU_{n}$. The collection $\{c U_{n}: n \in \mathbb{N}\}$ forms a neighbourhood base at $c \in Int(P)$. Hence there exists $N \in \mathbb{N}$ such that $cU_{n} \subset Int(P)$ for $n \geq N$.

Since $A$ is unital, it follows that $\int f(a)V_{a}^{*} \in A \rtimes_{red} P$. Thus , for $n \geq N$, $A \rtimes_{red} P$ contains  the product  $\displaystyle \Big(\int \alpha_{b}(x)\psi_{n}(b) V_{b}^{*} db \Big) \Big(\int f(a)V_{a}^{*}da \Big)=\int \int \alpha_{b}(x)\psi_{n}(b)f(a)V_{b}^{*}V_{a}^{*} da db$. Now observe that 
\begin{align*}
&\Big|\Big|\int \int \alpha_{b}(x)\psi_{n}(b)f(a)V_{b}^{*}V_{a}^{*} da db - \int \int \alpha_{c}(x)\psi_{n}(b)f(a)V_{b}^{*}V_{a}^{*} da db \Big|\Big|  \\
&\leq \displaystyle \sup_{b \in cU_{n}}(||\alpha_{b}(x)-\alpha_{c}(x)||) ||f||_{1}.
\end{align*}
Now the strong continuity of the action, together with the fact that $\{cU_{n}: n \in \mathbb{N}\}$ forms a neighbourhood base at $c$ with $\{cU_{n}\}$ decreasing implies that the difference \[\int \int \alpha_{b}(x)\psi_{n}(b)f(a)V_{b}^{*}V_{a}^{*} da db - \int \int \alpha_{c}(x)\psi_{n}(b)f(a)V_{b}^{*}V_{a}^{*} da db \rightarrow 0\]
in the norm.

Now calculate as follows to find that 
\begin{align*}
\int \int \alpha_{c}(x)\psi_{n}(b)f(a)V_{b}^{*}V_{a}^{*} da db & = \int \int \alpha_{c}(x)\psi_{n}(b)f(a)V_{ba}^{*} da db \\
& = \int \int \alpha_{c}(x) \psi_{n}(b)f(b^{-1}a)V_{a}^{*} da db \\
& = \int \alpha_{c}(x) \big (\int \psi_{n}(b)f(b^{-1}a)db \big) V_{a}^{*} da \\
& = \int \alpha_{c}(x) (\psi_{n}*f)(a)V_{a}^{*}da \\
& = \int \alpha_{c}(x)L_{c}(\phi_{n}*f)(a)V_{a}^{*} da \\
& \rightarrow \int \alpha_{c}(x)L_{c}(f)(a)V_{a}^{*} da.
\end{align*}
The last line follows from the fact that $(\phi_{n})$ is an approximate unit in $L^{1}(G)$ and left translations are continuous as bounded linear operators on $L^{1}(G)$. As a consequence, we have  $\int \int \alpha_{b}(x)\psi_{n}(b)f(a)V_{b}^{*}V_{a}^{*} da db \rightarrow \int \alpha_{c}(x)L_{c}(f)(a)V_{a}^{*}$. Since for $n \geq N$, $\int \int \alpha_{b}(x)\psi_{n}(b)f(a)V_{b}^{*}V_{a}^{*} da db \in A \rtimes_{red} P$, it follows that $\int \alpha_{c}(x)L_{c}(f)(a)V_{a}^{*}da \in A \rtimes_{red} P$. This proves our claim.

Now choose a sequence $c_{n} \in Int(P)$ such that $c_{n} \to e$. We leave it to the reader to prove that $\int \alpha_{c_{n}}(x)L_{c_{n}}(f)(a)V_{a}^{*}da \to \int x f(a)V_{a}^{*}da$. Thus $\int x f(a)V_{a}^{*}da \in A \rtimes_{red} P$. This completes the proof. \hfill $\Box$

Another way of stating the above lemma is that the Wiener-Hopf operators $W_{f} \in A \rtimes_{red} P$ if $f \in C_{c}(Int(P)^{-1},A)$. For let $f \in C_{c}(Int(P)^{-1},A)$ be given. Define $f_{1} \in C_{c}(G,A)$ by $f_{1}(g)=f(g^{-1})\Delta(g)^{-1}$. Then $f_{1} \in C_{c}(Int(P),A)$. Now note that
\begin{align*}
W_{f}&= \int f(g)W_{g}dg \\
    & = \int f(g^{-1})W_{g^{-1}}\Delta(g)^{-1} dg \\
    & = \displaystyle \int_{a \in Int(P)} f_{1}(a)V_{a}^{*}da \in A \rtimes_{red} P.
\end{align*}

\begin{ppsn}
\label{equality of Toeplitz}
The $*$-algebra generated by $\{\widetilde{f}: f \in C_{c}(Int(P),A)\}$ is dense in $\Gamma_{c}(\clg,r^{*}\cla)$ with respect to the $I$-norm. Consequently, $\mathfrak{T}(A,P,\alpha)=A \rtimes_{red} P$.
\end{ppsn}
\textit{Proof.} Let $\mathcal{D}_{+}$ be the $*$-algebra generated by $\{\widetilde{f}:f \in C_{c}(Int(P),A)\}$. Let $x \in A$ and $\phi,\psi \in C_{c}(Int(P))$ be given and let $f: = x \otimes \phi$. Define $\psi_{1} \in C_{c}(G)$ by $\psi_{1}(g)=\overline{\psi(g^{-1})}$. The map $C_{c}(Int(P)) \ni \psi \to \psi_{1} \in C_{c}(Int(P)^{-1})$ is a bijection. Now observe that $\widetilde{f}*\widetilde{\psi}^{*}= \widetilde{g}$ where $g=x \otimes (\phi * \psi_{1})$. Now the density of $\mathcal{D}_{+}$ follows from Lemma \ref{convolution_density} and Proposition \ref{density in inductive limit topology}.
 
 Clearly $A \rtimes_{red} P \subset \mathfrak{T}(A,P,\alpha)$.  Recall the representation $\Lambda_{X_{0}}$ considered in Theorem \ref{main_theorem}. The density of $\mathcal{D}_{+}$ in $\Gamma_{c}(\clg,r^{*}\cla)$ with respect to the $I$-norm implies that $*$-algebra generated by $\{\Lambda_{X_{0}}(\widetilde{f}): f \in C_{c}(Int(P),A)\}= \{ W_{f}: f \in C_{c}(Int(P)^{-1},A)\}$ is dense in $\mathfrak{T}(A,P,\alpha)$. But $W_{f} \in A \rtimes_{red} P$ if $f \in C_{c}(Int(P)^{-1},A)$. Hence $\mathfrak{T}(A,P,\alpha)\subset A \rtimes_{red} P$ and consequently $A \rtimes_{red} P= \mathfrak{T}(A,P,\alpha)$. This completes the proof. \hfill $\Box$

\section{The surjective case}

Throughout this section,   $A$ stands for a separable $C^{*}$-algebra, possibly non-unital, and   $\alpha:P \to End(A)$ is a strongly continuous surjective action i.e. for every $a \in P$, $\alpha_{a}$ is surjective. We do not demand that $P$ is left Ore. Let $\clg:=\Omega \rtimes P$ be the Wiener-Hopf groupoid. For $X \in \Omega$, let 
\[
I_{X}:= \{x \in A: P^{-1}a_{n} \to X,~a_{n} \in P \Rightarrow ||\alpha_{a_{n}}(x)|| \to 0 \}.
\]
Note that for $X \in \Omega$, $I_{X}$ is a closed two sided ideal of $A$. 

\begin{rmrk}
If $P^{-1}a=P^{-1}b$ for $a,b \in P$, then there exists $c \in P \cap P^{-1}$ such that $a=cb$. Hence for $x \in A$, $||\alpha_{a}(x)||=||\alpha_{b}(x)||$. This is because $\alpha_{c}:A \to A$ is an isomorphism if $c \in P \cap P^{-1}$.
\end{rmrk}

\begin{lmma}
\label{basic action}
Let $X \in \Omega$, $a \in P$ and $x \in A$ be given. Then $\alpha_{a}(x) \in I_{X}$ if and only if $x \in I_{Xa}$. 
\end{lmma}
\textit{Proof.} Let $X \in \Omega$, $a \in P$ and $x \in A$ be given. Suppose $x \in I_{Xa}$. Let $P^{-1}a_{n}$ be a sequence, with $a_{n} \in P$, in $\Omega$ converging to $X$. Then $P^{-1}a_{n}a \to Xa$. Thus by definition $\alpha_{a_{n}a}(x)=\alpha_{a_{n}}(\alpha_{a}(x)) \to 0$. Hence $\alpha_{a}(x) \in I_{X}$. This proves the "if" part.

We prove the "only if" part first for the case when $a \in Int(P)$. Thus assume $a \in Int(P)$ and $\alpha_{a}(x) \in I_{X}$. Suppose $P^{-1}a_{n} \to Xa$. Since $e \in X$, it follows that there exists $b_{n} \in P$ such that $b_{n}^{-1}a_{n} \to a$. But $a \in Int(P)$ and $Int(P)$ is open in $G$. Thus $c_{n}:=b_{n}^{-1}a_{n} \in Int(P)$ eventually. Since $P^{-1}a_{n} \to Xa$ and $b_{n}^{-1}a_{n} \to a$, it follows that $P^{-1}b_{n} \to X$. Hence $\alpha_{b_{n}}(\alpha_{a}(x)) \to 0$. Now observe that 
\begin{align*}
||\alpha_{a_{n}}(x)||&=||\alpha_{b_{n}c_{n}}(x)|| \\
                     & \leq ||\alpha_{b_{n}}(\alpha_{a}(x))|| + ||\alpha_{b_{n}}(\alpha_{c_{n}}(x)-\alpha_{a}(x))|| 
                     & \leq ||\alpha_{b_{n}}(\alpha_{a}(x))|| + ||\alpha_{c_{n}}(x)-\alpha_{a}(x)||.
\end{align*}
In the right hand side of the above inequality, we have just observed that the first term converges to zero. Since $c_{n} \to a$, the strong continuity of the action $\alpha$ implies that the second term also converges to zero. Hence $\alpha_{a_{n}}(x) \to 0$. This proves the "only if" part when $a \in Int(P)$.

Now let $a \in P$ and suppose $\alpha_{a}(x) \in I_{X}$. Choose $c \in Int(P)$. Since $\alpha_{c}$ is surjective, let $y \in A$ be such that $x=\alpha_{c}(y)$. Then $\alpha_{ac}(y) \in I_{X}$. But $ac \in Int(P)$. Thus by the "only if" part ( that is just proved), it follows that $y \in I_{Xac}$. Then the "if part" implies that $x=\alpha_{c}(y) \in I_{Xa}$. This completes the proof. \hfill $\Box$

\begin{lmma}
\label{basic facts regarding the ideals}
For $a \in P$, $I_{P^{-1}a}=Ker(\alpha_{a})$.
\end{lmma}
\textit{Proof.} Clearly the inclusion $I_{P^{-1}a} \subset Ker(\alpha_{a})$ holds. For the constant sequence $(P^{-1}a)_{n=1}^{\infty}$ tends to $ P^{-1}a$. Hence if $x \in I_{P^{-1}a}$ then $\alpha_{a}(x)=0$.
  
  We prove the reverse inclusion first for $a \in Int(P)$. Let $a \in Int(P)$ and $x \in Ker(\alpha_{a})$ be given. Consider a sequence $(P^{-1}a_{n})$, with $a_{n} \in P$, converging to $P^{-1}a$. Thus there exists a sequence $b_{n} \in P$ such that $c_{n}:=b_{n}^{-1}a_{n} \to a$. But $a \in Int(P)$. Hence $c_{n} \in Int(P)$ eventually. Write $a_{n}=b_{n}c_{n}$. Now $c_{n} \to a$ implies that $\alpha_{c_{n}}(x) \to \alpha_{a}(x)=0$. Now observe that for $n \geq 1$, $||\alpha_{a_{n}}(x)||=||\alpha_{b_{n}}\alpha_{c_{n}}(x)|| \leq ||\alpha_{c_{n}}(x)||$. It follows that $\alpha_{a_{n}}(x) \to 0$. That is $x \in I_{P^{-1}a}$.
  
  Now let $a \in P$ and $x \in Ker(\alpha_{a})$ be given. Choose $c \in Int(P)$ and write $x=\alpha_{c}(y)$ for some $y \in A$. Then $y \in Ker(\alpha_{ac})$ and $ac \in Int(P)$. By what we have proved just now, it follows that $y \in I_{P^{-1}ac}$. Now Lemma \ref{basic action} implies that $x=\alpha_{c}(y) \in I_{P^{-1}a}$. This proves the reverse inclusion and the proof is complete. \hfill $\Box$ 

\begin{rmrk}
Note that Lemma \ref{basic facts regarding the ideals} implies in particular that $I_{P^{-1}}=0$.
 \end{rmrk} 
For $X \in \Omega$, let $A_{X}:=A/I_{X}$ and let \[\cla:=\coprod_{X \in \Omega} A_{X}\] be the disjoint union and $p:\cla \to \Omega$ be the surjection sending $A_{X} \to X$. Thus $\cla$ is a bundle of $C^{*}$-algebras over $\Omega$. The next proposition is  crucial in endowing a topology on $\cla$ so that $\cla$ becomes an upper semi-continuous bundle.

We need a little lemma before proving it.

\begin{lmma}
\label{uniqueness of norms}
Let $B$ be a $C^{*}$-algebra with norm $||.||$. Let $||.||_{1}$ be norm on $B$ satisfying the $C^{*}$-identity $||x^{*}x||_{1}=||x||^{2}$. Suppose $||x||_{1} \leq ||x||$ for every $x \in B$ then $||x||=||x||_{1}$ for every $x \in B$.
\end{lmma}
\textit{Proof.} The point of the lemma is that we do not assume $(B,||.||_{1})$ is complete. Let $\overline{B}$ be the completion of $B$ w.r.t. the norm $||.||_{1}$. The hypothesis of the lemma implies that the inclusion $(B,||.||) \hookrightarrow (\overline{B},||.||_{1})$ is continuous, injective and is a $*$-homomorphism. Hence the image is closed. That is $(B,||.||_{1})$ is complete. Thus the two norms $||.||$ and $||.||_{1}$ on $B$ are equal.  This completes the proof. \hfill $\Box$

\begin{lmma}
We have the following. 
\begin{enumerate}
\item[(1)] Let $x \in A$ and $X \in \Omega$ be given. Then $||x+I_{X}||=\displaystyle \limsup_{P^{-1}a \rightarrow X}||\alpha_{a}(x)||$. \footnote{ Here $\displaystyle \limsup_{P^{-1}a \to X}||\alpha_{a}(x)||:= \inf_{X \in U \textrm{-open}}\big(\sup_{P^{-1}a \in U}||\alpha_{a}(x)||$~\big).} 
\item[(2)] For $x \in A$ the map $\Omega \ni X \to ||x+I_{X}|| $ is upper semi-continuous.
\end{enumerate}
\end{lmma}  
\textit{Proof.} Let $X \in \Omega$ be given. For $x \in A$, let \[\displaystyle ||x||_{1}:=\displaystyle \limsup_{P^{-1}a \rightarrow X}||\alpha_{a}(x)||.\] 

Note that $||.||_{1}$ is a semi-norm on $A$ satisfying the $C^{*}$-identity. Note that $||x||_{1} \leq ||x||$ for $x \in A$. We claim that  $I_{X}:= \{x \in A: ||x||_{1}=0\}$. Let $x \in A$ be such that $||x||_{1}=0$. Consider a sequence $(P^{-1}a_{n}) \to X$ and let $\epsilon >0$ be given. Since $||x||_{1}=0$, it follows that there exists an open set $U \subset \Omega$ containing $X$ such that $P^{-1}a \in U$ implies $||\alpha_{a}(x)|| \leq \epsilon$. Since $P^{-1}a_{n} \to X$, it follows that $P^{-1}a_{n} \in U$ eventually. Hence $||\alpha_{a_{n}}(x)|| \leq \epsilon$ eventually. This proves that $x \in I_{X}$.

Now let $x \in I_{X}$. Suppose $||x||_{1}>0$. Choose a countable base $\{U_{n}: n \in \mathbb{N}\}$ of open subsets of $\Omega$ containing $X$. We can assume that $U_{n}$ is decreasing. Since $||x||_{1}>0$, for every $n$, there exists $a_{n} \in P$ such that $P^{-1}a_{n} \in U_{n}$ and $||\alpha_{a_{n}}(x)|| > \frac{||x||_{1}}{2}$. Then $P^{-1}a_{n} \to X$ but $\alpha_{a_{n}}(x) \nrightarrow 0$. This contradicts the fact that $x \in I_{X}$. Hence $||x||_{1}=0$. This proves our claim. Consequently $||.||_{1}$ descends to a $C^{*}$-norm on $A/I_{X}$ which we again denote by $||.||_{1}$.

Note that for $x \in A$ and $y \in I_{X}$, $||x||_{1} =||x-y||_{1} \leq ||x-y||$. Thus $||x||_{1} \leq ||x+I_{X}||$ for $x \in A$ and consequently $||x+I_{X}||_{1} \leq ||x+I_{X}||$ for $x \in A$. Now by Lemma \ref{uniqueness of norms} applied to the $C^{*}$-algebra $A/I_{X}$ implies that $||x+I_{X}||_{1}=||x+I_{X}||$ for every $x$ in $A$. This proves $(1)$. 

Let $x \in A$ and $\epsilon>0$ be given. To show $(2)$, it is enough to show that  the set \[V:=\{X \in \Omega: ||x+I_{X}||<\epsilon\}\] is open in $\Omega$. Let $X_{0} \in V$. Then, by $(1)$, $\displaystyle \limsup_{P^{-1}a \to X_{0}}||\alpha_{a}(x)|| < \epsilon$. Choose $r>0$ such that $\displaystyle \limsup_{P^{-1}a \to X_{0}}||\alpha_{a}(x)|| < r<\epsilon$. Hence there exists an open set $U \subset \Omega$ containing $X_{0}$ such that $||\alpha_{a}(x)|| < r$ if $P^{-1}a \in U$. Hence $\displaystyle \sup_{P^{-1}a \in U}||\alpha_{a}(x)|| \leq  r < \epsilon$. Then for $X \in U$, $\displaystyle \limsup_{P^{-1}a \to X}||\alpha_{a}(x)|| \leq   r < \epsilon$. Hence $X_0\in U \subset V$. This proves that $V$ is open in $\Omega$. This proves $(2)$ and the proof is complete. \hfill $\Box$

 By Proposition C.25 in Appendix C of \cite{Dana-Williams}, it follows that there exists a unique topology $\tau$ on $\cla:=\displaystyle \coprod_{X \in \Omega}A_{X}$ such that $p:\cla \to \Omega$ is an upper semi-continuous bundle and the map $\Omega \in X \to x+I_{X}$ is a continuous section for every $x \in A$.  For an open subset $U \subset \Omega$, $x \in A$ and $\epsilon >0$, define 
 \[
 \mathcal{W}(U,x,\epsilon):= \{y+I_{X} \in \cla : X \in U \textrm{~and~} ||(y-x)+I_{X}|| < \epsilon\}. 
 \]
 Then the collection $\{\mathcal{W}(U,x,\epsilon)\}$ forms a basis for the topology $\tau$. In fact, since we have assumed that $A$ is separable, it follows that $\cla$ is second countable. 
 
Let $a \in P$ be given and $X \in \Omega$ be given. Now Lemma \ref{basic action} and the surjectivity of $\alpha_{a}$ implies that the map $A_{Xa} \ni x+I_{Xa} \to \alpha_{a}(x)+I_{X} \in A_{X}$ is an isomorphism. We denote this map by $\alpha_{(X,a)}$. Note that for $X \in \Omega$ and $a,b \in P$, $\alpha_{(X,a)}\alpha_{(Xa,b)}=\alpha_{(X,ab)}$.

 The Wiener-Hopf groupoid $\clg=\Omega \rtimes P$ acts on $\cla$ as follows. Let $(X,g) \in \clg$ be given and let $Y=Xg$. Write $g=ab^{-1}$ with $a,b \in P$. Define $\alpha_{(X,g)}:A_{Y} \to A_{X}$ by the formula $\displaystyle \alpha_{(X,g)}:=\alpha_{(X,a)}\alpha_{(Y,b)}^{-1}.$ First we need to verify that $\alpha_{(X,g)}$ is well-defined. 
 
 Suppose $g=a_{1}b_{1}^{-1}=a_{2}b_{2}^{-1}$. Then $a_{1}^{-1}a_{2}=b_{1}^{-1}b_{2}$. Choose $c,d \in P$ such that $a_{1}^{-1}a_{2}=cd^{-1}$. Then $a_{1}c=a_{2}d$. Similarly $b_{1}c=b_{2}d$. Also $Xa_{i}=Yb_{i}$ for $i=1,2$. 
 Note that 
 \begin{align*}
 \alpha_{(X,a_{1})}\alpha_{(Y,b_{1})}^{-1}&= \alpha_{(X,a_{1})}\alpha_{(Xa_{1},c)}\alpha_{(Xa_{1},c)}^{-1}\alpha_{(Y,b_{1})}^{-1} \\
 &= \alpha_{(X,a_{1})}\alpha_{(Xa_{1},c)}(\alpha_{(Y,b_{1})}\alpha_{(Yb_{1},c)})^{-1} \\
 &=\alpha_{(X,a_{1}c)}\alpha_{(Y,b_{1}c)}^{-1}
   \end{align*}
   Similarly $\alpha_{(X,a_{2})}\alpha_{(X,b_{2})}^{-1}=\alpha_{(X,a_{2}d)}\alpha_{(Y,b_{2}d)}^{-1}$. Since $a_{1}c=a_{2}d$ and $b_{1}c=b_{2}d$, it follows that $\alpha_{(X,g)}$ is well-defined.
   
   We leave it to the reader to verify that $\{\alpha_{(X,g)}: (X,g) \in \clg\}$ defines a left action of $\clg$ on $\cla$. In fact, the action is continuous, which is the content of the next proposition. 
   
 \begin{lmma}
 \label{continuity of the action}
 Let $(X_{n},a_{n},y_{n})$ be a sequence in $\Omega \times P \times A$ and let $(X,a,y) \in \Omega \times P \times A$. Set $x_{n}:=\alpha_{a_{n}}(y_{n})$ and $x:=\alpha_{a}(y)$. Suppose that $a_{n} \to a$. Then  $x_{n}+I_{X_{n}}\to x+I_{X}$ in $\cla$ if and only if  $y_{n}+I_{X_{n}a_{n}}\to y+I_{Xa}$ in $\cla$.
   \end{lmma}
  \textit{Proof.} Suppose $x_{n}+I_{X_{n}} \to x+I_{X}$ in $\cla$. Choose a basic open set $\mathcal{W}(U,z,\epsilon)$ of $\cla$ containing $y+I_{Xa}$.   Then $Xa \in U$ and $||(y-z)+I_{Xa}||< \epsilon$.   
  Since the map $\Omega \times P \ni (Z,b) \to Zb \in \Omega$ is continuous, there exists an open set $V \subset \Omega$ containing $X$ and an open subset $W \subset P$ (i.e. open in the subspace topology) containing $a$ such that $VW \subset U$. Choose $r>0$ such that $||(y-z)+I_{Xa}|| < r < \epsilon$. Since $a_{n} \to a$ and $X_{n} \to X$, there exists $N_{0} \in \mathbb{N}$ such that $n \geq N_{0}$ implies $X_{n} \in V$, $a_{n} \in W$ and $||\alpha_{a_{n}}(z)-\alpha_{a}(z)|| < \epsilon -r$. In particular, for $n \geq N_{0}$, $X_{n}a_{n} \in U$. 
  
  Since $\alpha_{(X,a)}:A_{Xa} \to A_{X}$ is an isomorphism, it follows that $||(x-\alpha_{a}(z))+I_{X}||< r$. Hence $x+I_{X} \in \mathcal{W}(V,\alpha_{a}(z),r)$. Since $x_{n}+I_{X_{n}} \to x+I_{X}$, it follows that there exists $N$ which we can assume to be larger than $N_{0}$ such that for $n \geq N$, $x_{n}+I_{X_{n}} \in \mathcal{W}(V,\alpha_{a}(z),r)$. 
  
 For $n \geq N$, $X_{n}a_{n} \in U$ and calculate as follows to find that
 \begin{align*}
 ||(y_{n}-z)+I_{X_{n}a_{n}}||& = ||(\alpha_{a_{n}}(y_{n})-\alpha_{a_{n}}(z))+I_{X_{n}}|| \\
                             & \leq ||(x_{n}-\alpha_{a}(z))+I_{X_{n}}|| + ||(\alpha_{a_{n}}(z)-\alpha_{a}(z))+I_{X_{n}}|| \\
                             & \leq ||(x_{n}-\alpha_{a}(z))+I_{X_{n}}|| + ||(\alpha_{a_{n}}(z)-\alpha_{a}(z))|| \\
                             & < r + \epsilon -r~~\big(\textrm{~for $x_{n}+I_{X_{n}} \in \mathcal{W}(V,\alpha_{a}(z),r)$}\big) \\
                             & < \epsilon.
                              \end{align*} 
 Thus for $n \geq N$, $y_{n}+I_{X_{n}a_{n}} \in \mathcal{W}(U,z,\epsilon)$. This proves the "only if" part.         
 
 Now suppose assume that $y_{n}+I_{X_{n}a_{n}} \to y+I_{Xa}$. First assume that $a \in Int(P)$. Then without loss of generality, we can assume that $a_{n} \in Int(P)$. Let $\mathcal{W}(U,z,\epsilon)$ be a basic open set containing $x+I_{X}$ and write $z=\alpha_{a}(\widetilde{z})$.  Now the assumption, $y_{n}+I_{X_{n}a_{n}} \to y+I_{Xa}$ implies that $X_{n}a_{n} \to Xa$. But $\Omega$ is compact and the right action of $P$ on $\Omega$ is injective. Hence it follows that $X_{n} \to X$. The fact that the open set $U$ contains $X$ implies that  there exists $N_{1}$ such that $n \geq N_{1}$ implies $X_{n} \in U$. 
 
 Since $x+I_{X} \in \mathcal{W}(U,z,\epsilon)$  there exists $r>0$ for which $||(x-z)+I_{X}|| < r < \epsilon$ i.e. $||\alpha_{a}(y)-\alpha_{a}(\widetilde{z})+I_{X}||<r$. Hence $||(y-\widetilde{z})+I_{Xa}|| < r$. Then $y+I_{Xa} \in \mathcal{W}(\Omega_{0},\widetilde{z},r)$ where $\Omega_{0}=\Omega Int(P)$. Recall that $\Omega_{0}$ is open in $\Omega$. Since $y_{n}+I_{X_{n}a_{n}} \to y+I_{Xa}$ and $a_{n} \to a$, choose $N \geq N_{1}$ such that for $n \geq N$, $y_{n}+I_{X_{n}a_{n}}\in \mathcal{W}(\Omega_{0},\widetilde{z},r)$ and $||\alpha_{a_{n}}(\widetilde{z})-\alpha_{a}(\widetilde{z})|| < \epsilon - r$.
 
 Now calculate as before to find that for $n \geq N$, 
 \begin{align*}
 ||(x_{n}-z)+I_{X_{n}}|| & \leq ||(\alpha_{a_{n}}(y_{n})-\alpha_{a_{n}}(\widetilde{z}))+I_{X_{n}}|| + ||(\alpha_{a_{n}}(\widetilde{z})-\alpha_{a}(\widetilde{z}))+I_{X_{n}}|| \\
 & \leq ||(y_{n}-\widetilde{z})+I_{X_{n}a_{n}}|| + ||\alpha_{a_{n}}(\widetilde{z})-\alpha_{a}(\widetilde{z})|| \\
 & < r + \epsilon -r ~~\big (~\textrm{for $y_{n}+I_{X_{n}a_{n}} \in \mathcal{W}(\Omega_{0},\widetilde{z},r)~\big)$} \\
 & < \epsilon.
 \end{align*} 
 This proves that for $n \geq N$, $x_{n}+I_{X_{n}} \in \mathcal{W}(U,z,\epsilon)$.      

Now suppose $a \in P$. Choose $c \in Int(P)$ and write $y_{n}=\alpha_{c}(\widetilde{y_{n}})$, $y=\alpha_{c}(\widetilde{y})$. Now the fact that $y_{n}+I_{X_{n}a_{n}} \to y+I_{Xa}$ and the "only if" part that we proved in the first half of the lemma implies that $\widetilde{y_{n}}+I_{X_{n}a_{n}c} \to \widetilde{y}+I_{Xac}$. But $ac \in Int(P)$ and $a_{n}c \to ac$. Thus by what we have proved just now ( for the "if part"), it follows that $\alpha_{a_{n}c}(\widetilde{y_{n}})+I_{X_{n}} \to \alpha_{ac}(\widetilde{y})+I_{X}$. Hence $x_{n}+I_{X_{n}} \to x+I_{X}$. This completes the proof. \hfill $\Box$

Just as in the injective case, we will denote both the action of $P$ on $A$ and the action $\{\alpha_{(X,g)}\}$ of $\clg$ on $\cla$ by $\alpha$.

\begin{ppsn}
The left action $\alpha:=\{\alpha_{(X,g)}:(X,g) \in \clg\}$ of $\clg$ on $\cla$ is continuous.
\end{ppsn}
\textit{Proof.} Suppose $(X_{n},g_{n}) \to (X,g) \in \clg$. Define $Y_{n}:=X_{n}g_{n}$ and $Y=Xg$. Let $y_{n}+I_{Y_{n}}$ be a sequence in $\cla$ converging to $y+I_{Y}$. To complete the proof, we need to show that the sequence $\displaystyle \alpha_{(X_{n},g_{n})}(y_{n}+I_{Y_{n}}) \to \alpha_{(X,g)}(y+I_{Y})$. Since $G=Int(P)Int(P)^{-1}$, there exists $a,b \in Int(P)$ such that $g=ab^{-1} \in Int(P)b^{-1}$ which is an open subset of $G$. Hence $g_{n} \in Int(P)b^{-1}$ eventually. Write $g_{n}=a_{n}b^{-1}$ for large $n$. Then $a_{n} \to a$. Choose $z_{n} \in A$ and $z \in A$ such that $y_{n}=\alpha_{b}(z_{n})$ and $y=\alpha_{b}(z)$.

Now Lemma \ref{continuity of the action} implies that $z_{n}+I_{Y_{n}b} \to z+I_{Yb}$. Since $Y_{n}b=X_{n}a_{n}$ and $Yb=Xa$, it follows that $z_{n}+I_{X_{n}a_{n}} \to z+I_{Xa}$. Now again by Lemma \ref{continuity of the action}, it follows that the sequence $\alpha_{a_n}(z_{n})+I_{X_{n}} \to \alpha_{a}(z)+I_{X}$. Thus by definition, it follows that 
$\alpha_{(X_{n},g_{n})}(y_{n}+I_{Y_{n}}) \to \alpha_{(X,g)}(y+I_{Y})$. This completes the proof. \hfill $\Box$

\textit{Notation:} For $x \in A$ and $g \in C_{c}(G)$, we denote the function $G \ni s \to xg(s) \in A$ by $x \otimes g$. 
 For $f \in C_{c}(G,A)$, let $\widetilde{f} \in \Gamma_{c}(\clg,r^{*}\cla)$ be defined by 
\[
\widetilde{f}(X,s):= f(s)+I_{X}
\]
for $(X,s) \in \clg$. From the definition of the topology on $\cla$, it follows that $\widetilde{f}$ is continuous when $f=x \otimes g$ for some $x \in A$ and $g \in C_{c}(G)$. Now the continuity of $\widetilde{f}$ for $f \in C_{c}(G,A)$ follows from the fact that the linear span of $\{x \otimes g: x \in A, g \in C_{c}(G)\}$ is dense in $C_{c}(G,A)$ in the inductive limit toplogy. 

Next we proceed towards proving a proposition analogous to Proposition \ref{density in  inductive limit topology}. 

\begin{lmma}
\label{self-adjointness}
Let $\mathcal{F}$ be the closure of the linear span of $\{\widetilde{f}: f \in C_{c}(G,A)\}$ in $\Gamma_{c}(\clg,r^{*}\cla)$ with respect to the $||~||_{I}$-norm. Then $\mathcal{F}$ is $*$-closed.
\end{lmma}
\textit{Proof.} It is enough to show that if $f=x \otimes g$ for some $x \in A$ and $g \in C_{c}(G)$ then $(\widetilde{f})^{*} \in \mathcal{F}$. Thus, let $x \in A$ and $g \in C_{c}(G)$ be given. Define $\psi \in \Gamma_{c}(\clg,r^{*}\cla)$ by the equation \[\psi(X,s)=(x^{*}+I_{X})\overline{g(s^{-1})}\]
 for $(X,s) \in \clg$. To complete the proof, it is enough to show that $\psi^{*} \in \mathcal{F}$.
 
 Let $K$ be the support of $g$. Since $G=Int(P)Int(P)^{-1}$ and $K$ is compact, it follows that there exists $b_{1},b_{2},\cdots,b_{n} \in Int(P)$ such that $K \subset \bigcup_{i}Int(P)b_{i}^{-1}$. Choose $b \in P$ such that $b \geq b_{i}$ for $i=1,2\cdots,n$. Then $b \in Int(P)$ and $K \subset Int(P)b^{-1}$. Choose an open set $U$ containing $K$ such that $\overline{U}$ is compact. Write $x=\alpha_{b}(y)$ for some $y \in A$. 
 
   Let $\epsilon >0$ be given. Since the action $\alpha:P \to End(A)$ is strongly continuous, for every $a \in Int(P)$, there exists an open set $V_{a} \subset Int(P)$ containing $a$ such that $||\alpha_{t}(y)-\alpha_{s}(y)|| \leq \epsilon$ for $t,s \in V_{a}$. Now $\bigcup_{a \in Int(P)}V_{a}b^{-1}$ contains $K$. Choose finitely many $V_{a}'s$, whose right translates by $b^{-1}$ cover $K$. Call them $V_{1},V_{2},\cdots,V_{n}$. By replacing $V_{i}$ by $V_{i} \cap Ub$, we can assume that $V_{i} \subset Ub$. For every $i$, choose $a_{i} \in V_{i}$. Then
   
    \begin{enumerate}
 \item[(a)]  $a_{i} \in V_{i}$ for $i=1,2,\cdots,n$,
 \item[(b)] if $t \in V_{i}$ then $||\alpha_{t}(y)-\alpha_{a_{i}}(y)|| \leq \epsilon$ for all $i$, and
 \item[(c)] the support $K \subset \bigcup_{i=1}^{n}V_{i}b^{-1}$ and $V_{i} \subset Ub$. 
 \end{enumerate} 
  Choose a partition of unity $\{\chi_{i}:i=1,2,\cdots,n\}$ such that $\chi_{i} \geq 0$,  $supp(\chi_{i}) \subset V_{i}b^{-1}$ and $\sum_{i=1}^{n} \chi_{i}=1$ on $K$. Let $\widetilde{\psi} \in \mathcal{F}$ be defined by 
 \[
 \widetilde{\psi}(X,s)=\sum_{i=1}^{n}(\alpha_{a_{i}}(y)+I_{X})g(s)\chi_{i}(s)
 \]
 for $(X,s) \in \clg$. Note that $supp(\widetilde{\psi}) \subset \Omega \times K$. For $(X,s) \in \clg$, calculate as follows to find that 
 \begin{align*}
 ||\psi^{*}(X,s)-\widetilde{\psi}(X,s)||&= || \sum_{i=1}\big(\alpha_{(X,s)}(x+I_{X.s})-(\alpha_{a_{i}(y)}+I_{X}))\chi_{i}(s)g(s)|| \\
 & = || \sum_{i=1}\big((\alpha_{(X,s)}(\alpha_{b}(y)+I_{X.s})-(\alpha_{a_{i}(y)}+I_{X}))\chi_{i}(s)g(s)|| \\
 &\displaystyle \leq  \sum_{i}||(\alpha_{sb}(y)-\alpha_{a_{i}}(y))+I_{X}||~ \chi_{i}(s)|g(s)|1_{V_{i}b^{-1}}(s)\\
 &\displaystyle \leq  \sum_{i}||(\alpha_{sb}(y)-\alpha_{a_{i}}(y))||1_{V_{i}}(sb)~ \chi_{i}(s)|g(s)|\\
 &   \leq \epsilon \sum_{i=1}\chi_{i}(s)|g(s)| \\
 & \leq \epsilon ||g||_{\infty}.
 \end{align*}
 This implies that there exists a sequence $\psi_{n} \in \mathcal{F}$ such that $\psi_{n} \to \psi^{*}$ in the inductive limit topology and hence in the $I$-norm. This completes the proof. \hfill $\Box$

\begin{rmrk}
\label{F_{+}}
Denote the closure of the linear span of $\{\widetilde{f}: f \in C_{c}(Int(P),A)\}$ and    the closure of the linear span of $\{\widetilde{f}: f \in C_{c}(Int(P)^{-1},A)\}$ w.r.t. the $I$-norm by $\mathcal{F}_{+}$ and $\mathcal{F}_{-}$ respectively. The proof of Lemma  \ref{self-adjointness} shows that 
$\mathcal{F}_{+}^{*} \subset \mathcal{F}_{-}$ and $\mathcal{F}_{-}^{*} \subset \mathcal{F}_{+}$. Hence $\mathcal{F}_{+}^{*}=\mathcal{F}_{-}$ and $\mathcal{F}_{-}^{*}=\mathcal{F}_{+}$. 
\end{rmrk}
 Let $\mathcal{D}$   be the $*$-subalgebra of $\Gamma_{c}(\clg,r^{*}\cla)$  generated by $\{\widetilde{f}: f \in C_{c}(G,A)\}$. Denote the closure of $\mathcal{D}$ in the $I$-norm by $\overline{\mathcal{D}}$. For $f \in C_{c}(G)$, let $\widetilde{f} \in C_{c}(\clg)$ be as defined in Theorem \ref{RW_theorem}. Let $\mathcal{D}_{0}$   be the $*$-subalgebra of $C_{c}(\clg)$  generated by $\{\widetilde{f}: f \in C_{c}(G)\}$. Denote the closure of $\mathcal{D}_{0}$ in the $I$-norm by $\overline{\mathcal{D}_{0}}$. Note that $\overline{\mathcal{D}_{0}}=C_{c}(\mathcal{G})$.
 
 We claim that $\overline{\mathcal{D}_{0}}*\overline{\mathcal{D}} \subset \overline{\mathcal{D}}$. Note that it is enough to show that for $f \in C_{c}(G)$ and $g \in C_{c}(G,A)$, $\widetilde{f}*\widetilde{g} \in \overline{\mathcal{D}}$ and $\widetilde{f}*(\widetilde{g})^{*} \in \overline{\mathcal{D}}$. In view of Lemma \ref{self-adjointness}, it is enough to show that $\widetilde{f}*\widetilde{g} \in \overline{\mathcal{D}}$ for $f \in C_{c}(G)$ and $g \in C_{c}(G,A)$. Again, we can assume $g$ is of the form $x \otimes \phi$ for some $\phi \in C_{c}(G)$. We can also assume $x$ is positive.
  
 Thus let  $f,g \in C_{c}(G)$ and  $x \in A$ be given. Suppose that $x$ is positive and write $x=y^*y$.  Let $\psi \in \Gamma_{c}(\clg,r^{*}\cla)$ be defined by \[\psi(X,s)=(x+I_{X})g(s)\]
  for $(X,s) \in \clg$. \textit{Claim: $\widetilde{f}*\psi \in \overline{\mathcal{D}}$.} 
 
 \textit{Proof of the claim.} Let $f_{1},g_{1} \in C_{c}(G,A)$ be defined by  
  $f_{1}(t)=yf(t^{-1})$, $g_{1}(t)=yg(t)$ for $t \in G$.  Then note that for $(X,s) \in \clg$, 
 \begin{align*}
 (\widetilde{f_{1}}^{*}*\widetilde{g_{1}})(X,s)&= \int \widetilde{f_{1}}^{*}(X,t)\alpha_{(X,t)}(y+I_{X.t})g(t^{-1}s)1_{X}(t^{-1})dt \\
                        & = \int \alpha_{(X,t)}(y^{*}+I_{X.t})\alpha_{(X,t)}(y+I_{X.t})f(t)g(t^{-1}s)1_{X}(t^{-1})dt \\
                        & = \int \alpha_{(X,t)}(x+I_{X.t})f(t)g(t^{-1}s)1_{X}(t^{-1})dt \\
                        & = (\widetilde{f}*\psi)(X,s)
  \end{align*}
  This proves that $\widetilde{f}*\psi \in \overline{\widetilde{D}}$. As a consequence, we have $C_{c}(\clg)*\overline{\mathcal{D}} \subset \overline{\mathcal{D}}$. 
 
 We now prove the analog of Proposition \ref{density in  inductive limit topology} in the surjective case.
 We will be brief in the proof as it is similar to that of Proposition \ref{density in inductive limit topology}.

For $\psi \in \Gamma_{c}(\clg, r^{*}\cla)$ and $a \in P$, let $R_{a}(\psi) \in \Gamma_{c}(\clg, r^{*}\cla)$ be defined by the formula
$
R_{a}(\psi)(X,s):=\alpha_{(X,a)}(\psi(Xa,a^{-1}s))
$ for $(X,s) \in \clg$. Then Lemma \ref{continuity of R_{a}} and Corollary \ref{continuity} holds in the surjective case as well. Since the proof is almost the same, word for word, we leave the verification to the reader. Recall also the following notation from the injective case.

\textit{Notation:} For $\phi \in C(\Omega)$ and $f \in C_{c}(G)$, let $\phi \otimes f \in C_{c}(\clg)$ be the function defined by  $(\phi \otimes f)(X,s)=\phi(X)f(s)$  for $(X,s) \in \clg$.  For $\phi \in C(\Omega)$ and $\psi \in \Gamma_{c}(\clg,r^{*}\cla)$,  let $\phi\cdot \psi \in \Gamma_{c}(\clg,r^{*}\cla)$ be the function defined by $(\phi\cdot\psi)(X,s)=\phi(X)\psi(X,s)$. We now prove the analog of Proposition \ref{density in  inductive limit topology} in the surjective case.
 We will be brief in the proof as it is similar to that of Proposition \ref{density in inductive limit topology}.

 \begin{ppsn}
 \label{density surjective case}
  The $*$-algebra generated by $\{\widetilde{f}: f \in C_{c}(G,A)\}$ is dense in $\Gamma_{c}(\clg,r^{*}\cla)$ in the topology induced by the norm $||~||_{I}$.
  \end{ppsn}
\textit{Proof.}  
Let $\mathcal{D}$ be the $*$-algebra generated by  $\{\widetilde{f}:f \in C_{c}(G,A)\}$ in $\Gamma_{c}(\clg,r^{*}\cla)$ and denote its closure w.r.t. the $I$-norm by $\overline{\mathcal{D}}$ . Note that sections of the form $\phi.\psi$ where $\phi \in C(\Omega)$ and $\psi$ is of the form $\psi(X,s)=(x+I_{X})f(s)$ where $x \in A$ and $f \in C_{c}(G)$ is dense in the $I$-norm. This follows from a partition of unity type argument and Proposition C.24 in Appendix C in \cite{Dana-Williams}.  Let $\phi \in C(\Omega)$, $x \in A$ and $f \in C_{c}(G)$ be given. Define $\psi \in \mathcal{D}$ by $\psi(X,s)=(x+I_{X})f(s)$. 
To complete the proof, it is enough to show that $\phi.\psi \in \overline{\mathcal{D}}$.

Using the fact that  $C_{c}(\clg)*\overline{\mathcal{D}} \subset \overline{\mathcal{D}}$ and arguing as in the injective case (Prop. \ref{density in inductive limit topology}), it follows that for $a \in Int(P)$, $\phi.R_{a}(\psi) \in \overline{\mathcal{D}}$. Now choose a sequence $a_{n} \in Int(P)$ such that $a_{n} \to e$. Then $\phi.R_{a_{n}}(\psi) \to \phi.\psi$ in the inductive limit topology and hence in the $I$-norm. Thus $\phi.\psi \in \overline{\mathcal{D}}$. This completes the proof. \hfill $\Box$ 
 
 The rest is similar to that of the injective case. Note that for $X_{0}=P^{-1}$, the fibre $A_{X_{0}}=A$. Arguing as in the injective case, we obtain that $\Lambda_{X_{0}}$ is a faithful representation of $\cla \rtimes_{red} \clg$ on $L^{2}(P,A)$ and for $f \in C_{c}(G,A)$, \[\Lambda_{X_{0}}(\widetilde{f})=W_{\widehat{f}}\] where $\widehat{f}(s)=f(s^{-1})\Delta(s)^{-\frac{1}{2}}$ and $W_{\widehat{f}}$ is the Wiener-Hopf operator with symbol $\widehat{f}$.

 Thus, together with Proposition \ref{density surjective case}, we obtain the following theorem
 
 \begin{thm}
 \label{surjective main_theorem}
 With the foregoing notation, the Toeplitz algebra $\mathfrak{T}(A,P,\alpha)$ is isomorphic to the reduced crossed product $\cla \rtimes_{red} \clg$.
 \end{thm}
 
 Just like in the injective case, we prove that $\mathfrak{T}(A,P,\alpha)=A \rtimes_{red} P$ under the hypotheses of the present section. We start with a little lemma. For $x \in A$ and $f \in C_{c}(Int(P),A)$, let $F_{x,f} \in C_{c}(Int(P),A)$ be defined by $F_{x,f}(a):=\alpha_{a}(x)f(a)$.
  \begin{lmma}
  \label{density of F_{x,f}}
  The span of $\{F_{x,f}: x \in A, f \in C_{c}(Int(P),A)\}$ is dense in $C_{c}(Int(P),A)$ where the latter is equipped with the inductive limit topology.
   \end{lmma}
 \textit{Proof.} It suffices to approximate functions of the type $x \otimes f$ where $x \in A$ and $f \in C_{c}(Int(P))$. Let $x \in A$ and $f \in C_{c}(Int(P))$ be given. Denote the support of $f$ by $K$. Choose an open set $U$ containing $K$ such that $U \subset Int(P)$ and $\overline{U}$ is compact. Let $\epsilon>0$ be given.
  For $c \in K$, choose $y_{c} \in A$ such that $\alpha_{c}(y_{c})=x$. Then for $c \in K$, there exists an open set $V_{c}$ containing $c$ such that $||\alpha_{b}(y_{c})-\alpha_{c}(y_{c})|| \leq \epsilon$ for $b \in V_{c}$. We can assume that $V_{c} \subset U$. Then $\displaystyle \bigcup_{c \in K} V_{c}$ contains $K$. Since $K$ is compact, there exists $c_{1},c_{2},\cdots,c_{n}$  such that $K \subset \bigcup_{i}V_{c_i}$. Choose a partition of unity $\{\phi_{i}:i=1,2\cdots,n\}$ such that $\sum_{i=1}^{n}\phi_{i}=1$ on $K$, $\phi_{i}\geq 0$ and $supp(\phi_{i}) \subset V_{c_i}$.
 We will denote $y_{c_i}$ by $y_{i}$ and $V_{c_i}$ by $V_{i}$. 
 
   Let $f_{i}:=\phi_{i}f$. Then $supp(F_{y_{i},f_{i}}) \subset U$. Now calculate as follows to conclude that for $a \in Int(P)$,
  \begin{align*}
  ||xf(a)-\sum_{i=1}^{n}F_{y_i,f_{i}}(a)||& = || \sum_{i=1}^{n}\phi_{i}(a)\alpha_{c_{i}}(y_{i})f(a)-\sum_{i=1}^{n}\alpha_{a}(y_{i})\phi_{i}(a)f(a) || \\
  & \leq \sum_{i=1}^{n} \phi_{i}(a)|f(a)|||\alpha_{c_{i}}(y_i)-\alpha_{a}(y_i)|| \\
  & \leq \sum_{i=1}^{n}\phi_{i}(a) |f(a)| \displaystyle \sup_{a \in V_{i}}||\alpha_{a}(y_i)-\alpha_{c_i}(y_i)|| \\
  & \leq \epsilon \sum_{i=1}^{n}\phi_i(a)|f(a)| \\
  & \leq \epsilon ||f||_{\infty}
    \end{align*}
   This completes the proof. \hfill $\Box$

Hence for $f \in C_{c}(Int(P),A)$, $\int f(a)V_{a}^{*}da \in A \rtimes_{red} P$. This is because for $x \in A$ and $f \in C_{c}(Int(P))$, $\int F_{x,f}(a)V_{a}^{*}da=\Big(\int x^{*}f(a)V_{a}da\Big)^{*} \in A \rtimes_{red} P$. In other words, the Wiener-Hopf operators $W_{f} \in A \rtimes_{red} P$ for $f \in C_{c}(Int(P)^{-1},A)$.

\begin{ppsn}
The $*$-algebra generated by $\{\widetilde{f}: f \in C_{c}(Int(P),A)\}$ is dense in $\Gamma_{c}(\clg,r^{*}\cla)$ in the topology induced by the $I$-norm. Thus $\mathfrak{T}(A,P,\alpha)=A \rtimes_{red} P$.
\end{ppsn}
\textit{Proof.} Let $\mathcal{D}_{+}$ denote the closure of the $*$-algebra generated by $\{\widetilde{f}: f \in C_{c}(Int(P),A)\}$ in $\Gamma_{c}(\clg,\cla)$ with respect to the $I$-norm.  Recall the notations $\mathcal{F}_{+}$, $\mathcal{F}_{-}$ from Remark \ref{F_{+}}. Then $\mathcal{F}_{+},\mathcal{F}_{-} \subset \mathcal{D}_{+}$. 
Let $x \in A$ be positive, $f \in C_{c}(Int(P))$ and $g \in C_{c}(Int(P)^{-1})$ be given. Write $x=y^{*}y$ for some $y \in A$.  By Lemma \ref{density of F_{x,f}}, it follows that there exists $z_{n} \in A$ and $f_{n} \in C_{c}(Int(P))$ such that $F_{z_n,f_n} \to y \otimes f$ in the inductive limit topology. Let $\phi_{n}:=y^{*}\otimes f_{n}$ and $\psi_{n}:=z_{n} \otimes g$. Note that $\phi_{n}$ and $\psi_{n} \in \mathcal{D}_{+}$.

Calculate as follows to find that for $(X,s) \in \clg$,
\begin{align*}
(\widetilde{\phi_{n}}*\widetilde{\psi_{n}})&(X,s) = \int \widetilde{\phi_{n}}(X,t)\alpha_{(X,t)}\big(\widetilde{\psi_{n}}(X.t,t^{-1}s\big)1_{X}(t^{-1})dt \\
&= \int_{t \in Int(P)} (y^{*}+I_{X})(\alpha_{t}(z_{n})+I_{X})f_{n}(t)g(t^{-1}s) \big[ \textrm{ For $supp(f_n) \subset Int(P) \subset X^{-1}$}\big] \\
& = \int (y^{*}+I_{X})\widetilde{F_{z_n,f_n}}(X,t)g(t^{-1}s)dt
\end{align*}
Now the fact that $F_{z_n,f_n}\to y \otimes f$ in the inductive limit topology   implies that $\widetilde{\phi_{n}}*\widetilde{\psi_{n}} \to \widetilde{\psi}$ in the inductive topology and consequently in the $I$-norm where $\psi:=x \otimes (f*g)$. Now Lemma \ref{convolution_density} and Proposition \ref{density surjective case} implies that $\mathcal{D}_{+}=\Gamma_{c}(\clg,r^{*}\cla)$.

Now the equality $\mathfrak{T}(A,P,\alpha)=A \rtimes_{red} P$ follows exactly as in Prop. \ref{equality of Toeplitz}. We leave the details to the reader. \hfill $\Box$
\begin{rmrk}
From now on, since we will only be considering injective (and unital) or surjective actions, we will not make the distinction between $\mathfrak{T}(A,P,\alpha)$ and $A \rtimes_{red} P$.
\end{rmrk}

\section{Morita equivalence} 
 
 The aim of this section is to show that $A \rtimes_{red} P$, in the two cases that we considered, is Morita-equivalent to a crossed product $C^{*}$-algebra. This follows from a construction, which is called the basic construction by Muhly and Williams in \cite{MW08}, and the fact that the Wiener-Hopf groupoid is equivalent to a transformation groupoid. First we recall the "basic construction" [Section 9.2, page 57] from \cite{MW08}. We need the notion of groupoid actions for which we refer the reader to \cite{MRW}. 
 
 Let $\clg$ and $\clh$ be locally compact Hausdorff second countable groupoids with Haar systems. Let $Z$ be a locally compact Hausdorff space which is a $(\clg,\clh)$ equivalence. This means that there exists two open surjective  maps $r:Z \to \mathcal{G}^{(0)}$ and $s:Z \to \clh^{(0)}$, a left $\mathcal{G}$-action i.e. a map 
 \[
 \mathcal{G}*Z: = \{(\gamma,z) \in \mathcal{G} \times Z: s(\gamma)=r(z)\} \to \gamma z \in Z
 \]
 and a right $\mathcal{H}$-action i.e. a map \[
 Z * \clh: = \{(z,\gamma) \in Z \times \clh : r(\gamma)=s(z)\} \to z\gamma \in Z \]
 such that the actions commute and  the "source" map $s:Z \to \clh^{(0)}$ and the "range" map $r:Z \to \clg^{(0)}$ descend to a homeomorphism from $\mathcal{G} \backslash Z$ to $\clh^{(0)}$ and a homeomorphism from $Z/\clh$ to $\clg^{(0)}$ respectively. Moreover both the $\mathcal{G}$ and $\clh$ actions are free and proper. Recall that a left $\clg$-action on a space $Z$ is said to be 
 \begin{enumerate}
 \item[(1)] free if $\gamma z=z$ for $\gamma \in \clg$ and $z \in Z$ then $\gamma \in \clg^{(0)}$.
 \item[(2)] proper if the map $\clg*Z \ni (\gamma,z) \to (\gamma z,z) \in Z \times Z$ is proper. 
 \end{enumerate}
 
 \begin{rmrk}
 We will denote the range and source maps of all the groupoids and the spaces which implement equivalences between them by $r$ and $s$ itself.
 \end{rmrk}

Let $Z$ be a $(\clg,\clh)$-equivalence and let $(\cla,\clg,\alpha)$ be a groupoid dynamical system. Define 
\[
r^{*}(\cla):=\{(z,a) \in Z \times \cla: a \in A_{r(z)} \}.
\]
The $r^{*}(\cla)$ is an upper semi-continuous bundle over $Z$ with $r^{*}(\cla) \ni (z,a) \to z \in Z$ being the open surjection. The groupoid $\clg$ acts on $r^{*}(\cla)$ as 
\[
\gamma.(z,a):=(\gamma z, \alpha_{\gamma}(a))
\]
if $s(\gamma)=r(z)$. 

Denote the quotient $\clg\backslash r^{*}(\cla)$ by $\cla_{Z}$. The map $\cla_{Z} \ni [(z,a)] \to s(z) \in \clh^{(0)}$ is well-defined and makes $\cla_{Z}$ into an upper semi-continuous bundle over $\clh^{(0)}$. Moreover the groupoid $\clh$ acts on $\cla_{Z}$ on the left by 
\[
\beta_{\eta}[(z,a)]:=[(z\eta^{-1},a)]
\]
if $s(\eta)=s(z)$. The triple $(\cla_{Z},\clh, \beta)$ is a groupoid dynamical system. The construction of the dynamical system $(\cla_{Z},\clh,\beta)$ is called the basic construction in \cite{MW08}. This construction had earlier appeared in \cite{MRWK} for continuous bundles.

If $\cla$ has enough sections( i.e. given $x \in \clg^{(0)}$ and $a \in A_{x}$, there exists $f \in \Gamma_{c}(\clg^{(0)},\cla)$ such that $f(x)=a$) then $\cla_{Z}$ has enough sections. Even though this fact is probably  known to experts the author is not able to find a proper reference and hence the following proposition is included in this paper. Observe that if $f:Z \to \cla$ is continuous, $f(z) \in \cla_{r(z)}$ and $f(\gamma z)=\alpha_{\gamma}(f(z))$ if $r(z)=s(\gamma)$ then the map \[\clh^{(0)} \ni s(z) \to [(z,f(z)] \in \cla_{Z}\] is a well-defined continuous section. Thus, with the foregoing notation, it suffices to prove the following proposition.

\begin{ppsn}
Let $z_{0} \in Z$ and $a \in \cla_{ r(z_0)}$ be given. Given $\epsilon>0$, there exists a continuous function $f:Z \to \cla$ such that 
\begin{enumerate}
\item[(1)] for $z \in Z$, $f(z) \in \cla_{r(z)}$,
\item[(2)] For $(\gamma,z) \in \clg*Z$, $f(\gamma.z)=\alpha_{\gamma}(f(z))$, and
\item[(3)] $||f(z_{0})-a|| \leq  \epsilon$.
\end{enumerate}
\end{ppsn}
\textit{Proof.} Let $A:=\{(\gamma,z): r(\gamma)=r(z)\} \subset \clg \times Z$ and $B:=\{(z,w):s(z)=s(w)\} \subset Z \times Z$. Note that the map $A \ni (\gamma,z) \to (\gamma^{-1}z,z) \in B$, which we denote by $\Phi$, is a homeomorphism. Choose a continuous section $\phi: \clg^{(0)} \to \cla$ such that $\phi(r(z_0))=a$. The map $\clg: \gamma \to \alpha_{\gamma}(\phi(s(\gamma))) \in \cla$ is continuous and its value for $\gamma=r(z_0)$ is $a$. Hence there exists an open set $U$ in $\clg$ containing $r(z_0)$ with compact closure such that $||\alpha_{\gamma}(\phi(s(\gamma)))-a|| < \epsilon$ if $\gamma \in U$ and $r(\gamma)=r(z_0)$.

Let $V \subset Z$ be open with compact closure containing $z_{0}$. Let $W$ be the image of $(U \times V) \cap A$ under $\Phi$. Then $W$ is an open subset of $B$ whose closure is compact and contains $(z_{0},z_{0})$. Choose $F \in C_{c}(B)$ such that $F \geq 0$, $F(z_{0},z_{0})>0$ and $supp(F) \subset W$.
Note that $M:=\int F(\gamma^{-1}z_{0},z_{0})d\lambda^{r(z_0)}(\gamma) >0$. By replacing $F$ by $\frac{F}{M}$, we can assume that the integral $\int F(\gamma^{-1}z_{0},z_0)d\lambda^{r(z_0)}(\gamma)=1$.

 Choose an extension $K \in C_{c}(Z \times Z)$ of $F$. Define $f:Z \to \cla$ by \[
 f(z):= \int K(\gamma^{-1}z,z_{0})\alpha_{\gamma}(\phi(s(\gamma)) d\lambda^{r(z)}(\gamma). \]
We leave it to the reader to verify that $f$ is continuous and $(1)$ and $(2)$ are satisfied. Condition $(2)$ is equivalent to the fact that if $(\gamma,z) \in A$, $f(\gamma^{-1}z)=\alpha_{\gamma}^{-1}f(z)$.  Now calculate as follows to find that 
\begin{align*}
||f(z_0)-a||=& || \int F(\gamma^{-1}z_{0},z_{0})\big(\alpha_{\gamma}\phi(s(\gamma))-a\big)d\lambda^{r(z_{0})}|| \\
             & \leq \int_{\gamma \in U} F(\gamma^{-1}z_{0},z_{0})||\alpha_{\gamma}(\phi(s(\gamma))-a)||d\lambda^{(r(z_0))} \\
             & \leq \epsilon  \int F(\gamma^{-1}z_{0},z_{0})d\lambda^{r(z_0)} \\
             & \leq \epsilon
\end{align*}
This completes the proof. \hfill $\Box$

We need the following theorem from \cite{SW13}[ See Corollary 19, Section 6]. The full version of it is proved in \cite{MW08}[ See Theorem 5.5 and Section 9.2]. 

\begin{ppsn}[\cite{SW13}]
\label{Sims-Williams}
Let $Z$ be a $(\clg,\clh)$-equivalence and $(\cla,\clg,\alpha)$ be a groupoid dynamical system. Then the reduced crossed products $\cla \rtimes_{red} \clg$ and $\cla_{Z} \rtimes_{red} \clh$ are Morita-equivalent.
\end{ppsn}
We recall a few more preliminaries about upper semi-continuous bundles. In particular, we need the following proposition to which the author claims no originality. We omit the proof as it follows easily from  Proposition C.20, Page 361 of \cite{Dana-Williams}. Let us start with a definition.
 
 \begin{dfn}
 Let $p:\cla \to X$ and $q:\clb \to Y$ be upper semi-continuous bundles. A pair $(\Phi,\phi)$ of maps where $\Phi:\cla \to \clb$ and $\phi:X \to Y$ is called a bundle map if 
 \begin{enumerate}
 \item [(1)] the maps $\Phi$ and $\phi$ are continuous,
 \item [(2)] for $x \in X$, $\Phi$ maps the fibre $A_{x}$ into $B_{\phi(x)}$ i.e. $q \circ \Phi = \phi \circ p$, and
 \item [(3)] for $x \in X$, the restriction $\Phi:A_{x} \to B_{\phi(x)}$ is a $*$-algebra homomorphism.
 \end{enumerate}  
  \end{dfn}
  
  \begin{ppsn}
  \label{fibrewise isomorphism}
  Let $p:\cla \to X$ and $q:\clb \to Y$ be upper semi-continuous bundles. Let $(\Phi,\phi):(\cla,X) \to (\clb,Y)$ be a bundle map. Suppose that $\phi$ is a homeomorphism and for every $x$, the restriction $\Phi:A_{x} \to B_{\phi(x)}$ is an isomorphism then $\Phi$ is a homeomorphism.
  \end{ppsn}
  
 If $p:\cla \to F$ is an upper semi-continuous bundle over $X$ and if $F \subset X$ then $p^{-1}(F)$ is an upper semi-continuous bundle over $F$ and is called the restriction of $\cla$ over $F$. We denote the restriction of $\cla$ over $F$ by $\cla|_{F}$. If $\cla$ has enough sections then $\cla|_{F}$ also has enough sections.

In what follows, let $\mathcal{G}:=\Omega \rtimes P$ be the Wiener-Hopf groupoid.  From Lemma \ref{equivalence}, the Wiener-Hopf groupoid $\mathcal{G}$ is equivalent to $\widetilde{\Omega} \rtimes G$ where $\displaystyle \widetilde{\Omega}:=\bigcup_{a \in P} \Omega a^{-1}$. The equivalence between $\mathcal{G}$ and $\widetilde{\Omega} \rtimes P$ is implemented by the space $Z:=\Omega \times G$ where the left action of $\mathcal{G}$ on $Z$ is given by
\[
(X,g).(Y,h):=(X,gh) \textrm{~if~}Xg=Y
\]
and the right action of $\widetilde{\Omega} \rtimes G$ on $Z$ is given by
\[
(Y,h).(Y^{'},g):=(Y,hg) \textrm{~if~}Yh=Y^{'}.
\]
 The range  $r:Z \to \Omega$ and the source $s:Z \to \widetilde{\Omega}$ maps are given by 
 \begin{align*}
 r(X,g)&:=X \\
 s(X,g)&:=X.g 
 \end{align*}
 for $(X,g) \in Z$. 
 
 Let $(p:\cla \to \Omega,\mathcal{G},\alpha)$ be a groupoid dynamical system  and let $(\cla_{Z},\widetilde{\Omega} \rtimes G, \beta)$ be the dynamical system obtained from the basic construction where $Z:=\Omega \times G$ is the space implementing the equivalence between $\clg$ and $\widetilde{\Omega} \rtimes G$. Let $\widetilde{\beta}$ be the action of $G$ on $\Gamma_{0}(\widetilde{\Omega},\cla_{Z})$ given by the formula: For $f \in \Gamma_{0}(\widetilde{\Omega},\cla_{Z})$, $g \in G$ and $X \in \widetilde{\Omega}$, \[\widetilde{\beta}_{g}f(X):=\beta_{(X,g)}f(X.g).\]
 For $a \in P$, let the map $\Omega \ni X \to Xa^{-1} \in \widetilde{\Omega}$ be denoted by $j_{a}$ and let $\Phi_{a}:\cla \to \cla_{Z}$ be defined by $\Phi_{a}(x):=[((X,a^{-1}),x)]$ where $X:=p(x)$. Note that $(\Phi_{a},j_{a}):(\cla,\Omega) \to (\cla_{Z}|_{\Omega a^{-1}},\Omega a^{-1})$ is a bundle map. 
  
\begin{ppsn}
\label{Morita-eq}
With the foregoing notations, we have the following.
\begin{enumerate}
\item[(1)] For every $a \in \cla$, the bundle map $(\Phi_{a},j_{a}):(\cla,\Omega) \to (\cla_{Z}|_{\Omega a^{-1}},\Omega a^{-1})$ is a homeomorphism.
\item[(2)] The reduced crossed product $\cla \rtimes_{red} \clg$ is Morita-equivalent to the crossed product $\Gamma_{0}(\widetilde{\Omega},\cla_{Z}) \rtimes_{red} G$.
\end{enumerate}
\end{ppsn} 
\textit{Proof.} Let $a \in P$ be given. Clearly $\Phi_{a}$ is continuous. It is enough to show that $\Phi_{a}$ is an isomorphism fibrewise. Let $X \in \Omega$ be given. Since the the left action of $\clg$ on $Z$ is free, it follows that $\Phi_{a}$ is $1$-$1$. Now let $[(Y,g,y)] \in \cla_{Z}|_{\Omega a^{-1}}$ be such that $Yg=Xa^{-1}$. Then $(X,a^{-1}g^{-1}) \in \clg$ as $Y \in \Omega$. Now $(X,a^{-1}g^{-1}).(Y,g,y)=(X,a^{-1},\alpha_{(X,a^{-1}g^{-1})}(y))$. Thus $[(Y,g,y)]:=\Phi_{a}(x)$ where $x:=\alpha_{(X,a^{-1}g^{-1})}(y)$. By Proposition \ref{fibrewise isomorphism}, it follows that $\Phi_{a}$ is a homeomorphism. Statement $(2)$ follows from Proposition \ref{Sims-Williams} and Remark \ref{some remarks}. This completes the proof. \hfill $\Box$

We apply the above to the reduced crossed product  $A \rtimes_{red} P$. Let $A$ be a $C^{*}$-algebra and $\alpha:P \to End(A)$ be a strongly continuous semigroup homomorphism. Assume one of the following holds.
\begin{enumerate}
\item[(1)] The algebra $A$ is unital, $P$ is left Ore and the action $\alpha$ is injective, or
\item[(2)] The action $\alpha$ is surjective. 
\end{enumerate} 
In both cases, we denote the groupoid dynamical system constructed in Section 5 and 6 by $(\cla,\clg,\alpha)$ i.e. $(\cla,\clg,\alpha)$ stands for the dynamical system constructed in Section 5 (Section 6) if $\alpha$ is injective (surjective). Denote by $(\cla_{Z},\widetilde{\Omega} \rtimes G,\beta)$ the groupoid dynamical system constructed by applying the basic construction to $(\cla,\clg,\alpha)$.

\begin{crlre}
\label{Mor-eq}
With the above assumptions and notations, the reduced crossed product $A \rtimes_{red} P$ is Morita-equivalent to $\Gamma_{0}(\widetilde{\Omega},\cla_{Z})\rtimes_{red} G$.
\end{crlre}

We finish this section, by comparing the results obtained in this paper with those in \cite{Li13}. Let us recall the setup considered  in \cite{Li13}. 
For the rest of this section, let  $G$ be a discrete group and $Q \subset G$ be a subsemigroup containing $e$. Let $A$ be a $C^{*}$-algebra and $\alpha:G \to Aut(A)$ be an action. Consider the Hilbert $A$-module $F:=A \otimes \ell^{2}(Q)$. For $x \in A$, let $\widetilde{\pi}(x) \in \mathcal{L}_{A}(F)$ be the operator defined by the equation
\[
\widetilde{\pi}(x)(y \otimes \delta_{b}):=\alpha_{b}^{-1}(x)y \otimes \delta_{b}.
\]
For $a \in Q$, let $\widetilde{V}_{a}$ be the operator on $\mathcal{L}_{A}(F)$ defined by the equation
\[
\widetilde{V}_{a}(y \otimes \delta_{b})=y \otimes \delta_{ab}.
\]
Here $\{\delta_{b}: b \in Q\}$ denotes the standard orthonormal basis of $\ell^{2}(Q)$.  Li defines the reduced crossed product (Defn.3.3, \cite{Li13}), let us denote it by $A \rtimes_{\ell} Q$, as the $C^{*}$-subalgebra of $\mathcal{L}_{A}(F)$ generated by $\{\widetilde{\pi}(x)\widetilde{V}_{a}: x \in A, a \in Q\}$. 

\begin{rmrk}
The subscript $\ell$ in the notation $A \rtimes_{\ell} Q$ stands for the fact that Li considers the "left regular representation" of the semigroup $Q$ whereas we consider the right regular one.  
\end{rmrk}
Set $P=Q^{-1}$. For $x \in A$ and $a \in P$, recall from Section 4, the operators $\pi(x)$ and $V_{a}$,   on the Hilbert $A$-module $E=A \otimes \ell^{2}(P)$, defined by the equations
\begin{align*}
\pi(x)(y \otimes \delta_{b})&=\alpha_{b}(x)y \otimes \delta_{b} \\
V_{a}(y \otimes \delta_{b})&=y \otimes \delta_{ba}. 
\end{align*}
Let $U:F \to E$ be the unitary defined by $U(y \otimes \delta_{b})=y \otimes \delta_{b^{-1}}$. Note that 
\begin{align*}
U^{*}\pi(x)U&=\widetilde{\pi}(x) \\
U^{*}V_{a}U&=\widetilde{V}_{a^{-1}}
\end{align*}
for $x \in A$ and $a \in P$. Thus $A \rtimes_{\ell} Q$ is isomorphic to the reduced crossed product $A \rtimes_{red} P$.

\begin{rmrk}
 In \cite{Li13}, $A \rtimes_{\ell} Q$ is realised as a groupoid crossed product (Thm. 5.24, \cite{Li13}). We convince ourselves that the groupoid dynamical system constructed in \cite{Li13} coincides with 
the one that we have constructed in Section 6.   We should note that Li also considers semigroups which are not of Ore type.
\end{rmrk}
From now on assume that $Q$ is left Ore or equivalently $P$ is right Ore  i.e. $PP^{-1}=G$.  

Let $\clg:=\Omega \rtimes P$ be the Wiener-Hopf groupoid and let $\widetilde{\Omega}=\bigcup_{a \in P} \Omega a^{-1}$.  Note by Lemma \ref{equivalence} that $\Omega$ is a clopen subset of $\widetilde{\Omega}$. 
Let $(\cla,\clg,\alpha)$ be the groupoid dynamical system constructed in Section 6. Note that since $\alpha_{a}$ is an automorphism for every $a \in P$, the bundle $\cla$ is isomorphic to the trivial bundle $A \times \Omega$. 
The action $\alpha$ of the groupoid $\clg$ on $A \times \Omega$ is given by $\alpha_{(X,g)}(x,Xg)=(\alpha_{g}(x),X)$. 

 Let $(\cla_{Z},\widetilde{\Omega} \rtimes G, \beta)$ be the bundle obtained by applying the basic 
construction to $(A \times \Omega, \clg, \alpha)$ where $Z=\Omega \times G$ is the space implementing the equivalence between $\clg$ and $\widetilde{\Omega} \rtimes G$. Consider the trivial bundle $A \times \widetilde{\Omega}$ over $\widetilde{\Omega}$.
The transformation groupoid $\widetilde{\Omega} \rtimes G$ acts on $A \times \widetilde{\Omega}$, let us denote the action by $\widetilde{\alpha}$, as follows: $\widetilde{\alpha}_{(X,g)}(x,Xg)=(\alpha_{g}(x),X)$. 
We leave it to the reader to verify that the map $\mathcal{A}_{Z} \ni [((X,g),x)] \to (Xg,\alpha_{g}^{-1}(x)) \in A \times \widetilde{\Omega}$ is an isomorphism and is $\widetilde{\Omega} \rtimes G$-equivariant. Thus we identify $(\mathcal{A}_{Z},\widetilde{\Omega} \rtimes G, \beta)$ with
$(A \times \widetilde{\Omega}, \widetilde{\Omega} \rtimes G, \widetilde{\alpha})$. 

Next we recall the groupoid dynamical system constructed by Li in \cite{Li13}. Denote the $C^{*}$-algebra of bounded complex valued functions on $G$, with the supremum norm, by $\ell^{\infty}(G)$.  For $g \in G$, let $\gamma_{g}: \ell^{\infty}(G) \to \ell^{\infty}(G)$ be defined by 
$\gamma_{g}(f)(x)=f(g^{-1}x)$. Then $\gamma$ is a left action of $G$ on $\ell^{\infty}(G)$. For $A \subset G$, let $1_{A} \in \ell^{\infty}(G)$ be the characteristic function associated to $A$.  

Let $\mathcal{D}$ be the  $\gamma$-invariant $C^{*}$-subalgebra of $\ell^{\infty}(G)$ generated by $1_{Q}$.
Denote the spectrum of $\mathcal{D}$ by $\widetilde{\Omega}^{'}$. The group $G$ acts on the right on $\widetilde{\Omega}^{'}$ as follows: for $g \in G$ and $\chi \in \widetilde{\Omega}^{'}$, $\chi.g=\chi \circ \gamma_{g}$. Let $N:=\{\chi \in \widetilde{\Omega}^{'}: \chi(1_{Q})=1\}$. 
Let $\mathcal{G}^{'}$ be the restriction of the transformation groupoid $\widetilde{\Omega}^{'} \rtimes G$ onto $N$. 
The transformation groupoid $\widetilde{\Omega}^{'} \rtimes G$ acts on the trivial bundle 
$A \times \widetilde{\Omega}^{'}$, let us denote the action by $\alpha^{'}$, by the formula : $\alpha^{'}_{(\chi,g)}(x,\chi.g)=(\alpha_{g}(x),\chi)$. The restricted groupoid $\clg^{'}$ acts on $A \times N$ via $\alpha^{'}$ and we denote the action of $\clg^{'}$ on $A \times N$ by $\alpha^{'}$ itself. In \cite{Li13} (Thm. 5.24.), it is shown that $A \rtimes_{\ell} Q$ is isomorphic to the reduced crossed product $(A \times N) \rtimes_{red} \clg^{'}$. 

Thus to convince ourselves that the groupoid dynamical system obtained in \cite{Li13} and by us are the same, it is enough to prove that there exists a $G$-equivariant homeomorphism from $\widetilde{\Omega}$ onto $\widetilde{\Omega}^{'}$ mapping $\Omega$ onto $N$. 
Recall the left action of $G$ on $C_{0}(\widetilde{\Omega})$. For $g \in G$ and $f \in C_{0}(\widetilde{\Omega})$, let $R_{g}(f) \in C_{0}(\widetilde{\Omega})$ be defined by $R_{g}(f)(X)=f(Xg)$.  Also note by Lemma \ref{equivalence} that $\Omega$ and its translates 
are clopen subsets of $\widetilde{\Omega}$. 

Let $\epsilon:C_{0}(\widetilde{\Omega}) \to \ell^{\infty}(G)$ be defined by $\epsilon(f)(g)=f(P^{-1}g^{-1})$. Clearly $\epsilon$ is $G$-equivariant. Since $\{P^{-1}g: g \in G\}$ is dense in $\widetilde{\Omega}$, it follows that
$\epsilon$ is injective.  By $(4)$ of Lemma \ref{equivalence}, it follows that $\epsilon(1_{\Omega})=1_{Q}$. Note again by Lemma \ref{equivalence} that $1_{\Omega g}(A)=1_{A}(g)$ for $A \in \widetilde{\Omega}$ and
$g \in G$. By the Stone-Weierstrass theorem, it follows that the $C^{*}$-subalgebra generated by $\{1_{\Omega g}: g \in G\}$ is dense in $C_{0}(\widetilde{\Omega})$. Consequently, the image of $\epsilon$ is $\mathcal{D}$. 
Thus the map $\epsilon:C_{0}(\widetilde{\Omega})\to \mathcal{D} \cong C_{0}(\widetilde{\Omega}^{'})$ is a $G$-equivariant isomorphism such that $\epsilon(1_{\Omega})=1_{Q}$. Denote the $G$-equivariant homeomorphism from $\widetilde{\Omega}$ to $\widetilde{\Omega}^{'}$ 
induced by the map $\epsilon$ by $\widetilde{\epsilon}$. Then clearly $\widetilde{\epsilon}$ maps $\Omega$ onto $N$. 

As a result, we conclude that our construction (when the action of the semigroup $P$ on the $C^{*}$-algebra $A$ is given by an action of the group $G$) coincides with Li's construction in \cite{Li13}.

\section{K-group computations}
We say that a locally compact group $H$ has the Connes-Thom isomorphism property with period $n$ if whenever $A$ is a $C^{*}$-algebra and $\alpha:H \to Aut(A)$ is a strongly continuous action then $K_{i+n}(A \rtimes_{red,\alpha} H) \cong K_{i}(A)$.

For the rest of this paper, we assume that the semigroup $P$ is \textbf{connected} and the group $G$ has the \textbf{Connes-Thom isomorphism} property with period $n$. The following are examples of such pairs.
\begin{enumerate}
\item[(1)] Let $G=\mathbb{R}^{n}$ and $P$ be a closed convex cone such that $P-P=\mathbb{R}^{n}$. 
\item[(2)]  Let \[G:=\Big \{\begin{pmatrix}
                                                   a & b \\
                                                   0 & 1
                                                  \end{pmatrix}: a > 0 , b \in \mathbb{R} \Big\}\] be the $ax+b$-group. Note that $G$ is isomorphic to   $\mathbb{R} \rtimes (0,\infty)$. Hence $G$ has the Connes-Thom property with period $2$.
                  Let \[P=\Big \{\begin{pmatrix}
                                                   a & b \\
                                                   0 & 1
                                                  \end{pmatrix}: a  \geq 1, b  \geq 0  \Big\}.\] For an explicit model for the order compactification of $P$, we refer the reader to \cite{Jean_Sundar}.
 \item[(3)] Let $G$ be the Heisenberg group $\mathbb{R}^{2} \rtimes \mathbb{R}$ where $\mathbb{R}$ acts on $\mathbb{R}^{2}$ by $a.(x,y)=(x,y+ax)$ and let $P:=\{(x,y,a) \in G: x,y,a \geq 0\}$. Note that $G$ has the Connes-Thom isomorphism property with period 3. The order compactification of $P$ is explicitly described in \cite{Nica90}
\end{enumerate}

 Let $\clg:=\Omega \rtimes P$ be the Wiener-Hopf groupoid and let $(\cla,\clg,\alpha)$ be a groupoid dynamical system. 
  
  Denote by $(\cla_{Z},\widetilde{\Omega} \rtimes G,\beta)$ the groupoid dynamical system constructed by applying the basic construction to $(\cla,\clg,\alpha)$ where $Z:=\Omega \times G$ and $\widetilde{\Omega}:=\bigcup_{a \in P} \Omega a^{-1}=\bigcup_{a \in Int(P)} \Omega_{0}a^{-1}$ where $\Omega_{0}=\Omega Int(P)$. 
 
 If $\widetilde{\mathcal{A}}$ is an upper semicontinuous bundle over a locally compact space $Y$ and if $U \subset Y$ is open, we denote the space of continuous sections of the bundle $\widetilde{\mathcal{A}}\to Y$ vanishing outside $U$ by $\Gamma_{0}(U,\widetilde{\mathcal{A}})$. Note that $\Gamma_{0}(U,\widetilde{\mathcal{A}})$ is the closure of $\{s \in \Gamma_{c}(Y,\widetilde{\mathcal{A}}): supp(s) \subset U\}$.
 
 \begin{thm}
 \label{K-theory}
 With the foregoing notations and assumptions, we have that the $K$-group $K_{i}(\cla \rtimes_{red} \clg)$  is isomorphic to $K_{i+n}(\Gamma_{0}(\Omega_{0},\cla))$ for $i=0,1$.
\end{thm}
\textit{Proof.} By Proposition \ref{Morita-eq} and from the fact that $G$ satisfies the Connes-Thom isomorphism property with period $n$, it follows that $K_{i}(\cla \rtimes_{red} \clg) \cong K_{i+n}(\Gamma_{0}(\widetilde{\Omega},\cla_{Z}))$. Recall that the action $\widetilde{\beta}$ of $G$  on $\Gamma_{0}(\widetilde{\Omega},\cla_{Z})$ is given by: for $g \in G$ and $f \in \Gamma_{0}(\widetilde{\Omega},\cla_{Z})$, $\widetilde{\beta}_{g}(f)(X):=\beta_{(X,g)}f(Xg)$. By Proposition \ref{Morita-eq}, we can view $\Gamma_{0}(\Omega_{0},\cla)$ as a $*$-subalgebra of $\Gamma_{0}(\widetilde{\Omega},\cla_{Z})$ via the bundle map $(\Phi_{e},j_{e})$. Moreover for $a \in P$, $\widetilde{\beta}_{a^{-1}}$ leaves $\Gamma_{0}(\Omega_{0},\cla)$ invariant. 

Note that $(P,\leq)$ is directed and if $a \leq b$ with $a,b \in P$ then $\Omega_{0}a^{-1} \subset \Omega_{0}b^{-1}$ Also $\widetilde{\Omega}:=\bigcup_{a \in P}\Omega_{0}a^{-1}$. Thus $\Gamma_{0}(\widetilde{\Omega},\cla_{Z})$ is the inductive limit of $\Gamma_{0}(\Omega_{0}a^{-1},\cla_{Z}|_{\Omega_{0}a^{-1}})$ where if $a \leq b$ with $a,b \in P$, the connecting map from $\Gamma_{0}(\Omega_{0}a^{-1},\cla_{Z}|_{\Omega_{0}a^{-1}}) \to \Gamma_{0}(\Omega_{0}b^{-1},\cla_{Z}|_{\Omega_{0}b^{-1}})$ is given by the inclusion map and we denote it by $i_{b,a}$.
 
By Proposition \ref{Morita-eq}, for $a \in P$, we can identify $\Gamma_{0}(\Omega_{0}a^{-1},\cla_{Z}|_{\Omega_{0}a^{-1}})$ with $\Gamma_{0}(\Omega_{0},\cla)$. After this identification, the connecting map $i_{b,a}$ is just the restriction of $\widetilde{\beta}_{b^{-1}a}$ to $\Gamma_{0}(\Omega_{0},\cla)$. Thus to complete the proof, it is enough to show that for $a \in P$, $\widetilde{\beta}_{a^{-1}}$ induces the identity map at the level of $K$-theory of $\Gamma_{0}(\Omega_{0},\cla)$. 

Let $a \in P$ be given and let $\phi_{t}(a)$ be a path in $P$ connecting $a$ to the identity element $e$. Then the restriction of $\widetilde{\beta}^{-1}_{\phi_{t}(a)}$ to $\Gamma_{0}(\Omega_{0},\cla)$  is a homotopy connecting $\widetilde{\beta}_{a^{-1}}$ to the identiy map. Now the theorem follows from the fact that $K$-theory preserves inductive limits. This completes the proof. \hfill $\Box$

Let us denote the set $\Omega\backslash \Omega_{0}$ by $\partial \Omega$. Let $\tau:P \to \Omega$ be the map given by $\tau(a)=P^{-1}a$.
We say that the semigroup $P$ satisfies condition \textbf{(H)} if there exists a homotopy of continuous maps $(\phi_{t})_{t \in [0,1]}:\Omega \to \Omega$ such that 
\begin{enumerate}
\item[(1)] For $t \in [0,1]$, $\phi_{t}$ leaves $\partial \Omega$ invariant,
\item[(2)] For $X \in \Omega$ and $t \in (0,1]$, $\phi_{t}(X) \in \tau(P)$ for $t \in (0,1]$ and $\phi_{t}(X) \subset  X$,
\item[(3)] $\phi_{0}=id_{\Omega}$, and $\phi_{1}(\Omega) \subset \partial \Omega$.
\end{enumerate}

\begin{ppsn}
\label{K-theory-injectivity}
Let $A$ be a unital $C^{*}$-algebra and $\alpha:P \to End(A)$ be an unital injective action. If the semigroup $P$ satisfies condition \textbf{(H)} then $K_{i}(A \rtimes_{red,\alpha} P)=0$ for $i=0,1$.
\end{ppsn}
\textit{Proof.} Let $(\cla,\clg,\alpha)$ be the groupoid dynamical system constructed in Section 5. Note that if $X \subset Y$ with $X,Y \in \Omega$ then $A_{X} \subset A_{Y}$. Also observe that \[\Gamma(\Omega,\cla):=\{f:\Omega \to B: \textrm{$f$ is continuous and }  f(X) \in A_{X} \textrm{ for $X \in \Omega$}\}.\] Then $\Gamma_{0}(\Omega_{0},\cla)$ is the subalgebra of sections which vanish on $\partial \Omega$. By Theorem \ref{K-theory}, to complete the proof, it suffices to prove that $\Gamma_{0}(\Omega_{0},\cla)$ has trivial $K$-theory.

Let $(\phi_{t})_{[0,1]}$ be a homotopy given by Condition \textbf{(H)}. For $t \in [0,1]$, define the homomorphism $\pi_{t}:\Gamma(\Omega,\cla) \to \Gamma(\Omega,\cla)$ by the equation $\pi_{t}(f)(X):=f(\phi_{t}(X))$ for $f \in \Gamma(\Omega,\cla)$.  Since $\phi_{t}(X) \subset X$ for $X \in \Omega$ and $t \in [0,1]$, it follows that $\pi_{t}(f) \in \Gamma(\Omega,\cla)$ for $f \in \Gamma(\Omega,\cla)$. Moreover $(1)$ of Condition \textbf{(H)} implies that $\pi_{t}$ leaves $\Gamma_{0}(\Omega_{0},\cla)$ invariant. Now, by $(3)$ of Condition \textbf{(H)}, $(\pi_{t})_{t \in [0,1]}$ restriced to $\Gamma_{0}(\Omega_{0},\cla)$ is a homotopy of $*$-algebra homomorphisms connecting the identity map with the zero map. Hence $\Gamma_{0}(\Omega_{0},\cla)$ has trivial $K$-theory. This completes the proof. \hfill $\Box$

Now let $A$ be a $C^{*}$-algebra and $\alpha:P \to End(A)$ be a surjective action. Let $(\cla,\clg,\alpha)$ be the dynamical system constructed in Section 6. Recall that for $X \in \Omega$ the fibre of $\cla$ at $X$ is given by the quotient $A_{X}:=A/I_{X}$ where the ideal $I_{X}$ is given by \[I_{X}:=\{x \in A: P^{-1}a_{n} \to X \Rightarrow \alpha_{a_{n}}(x) \to 0\}.\]

We claim that for $X \in \Omega$ and $a \in X$, $I_{P^{-1}a} \subset I_{X}$. Suppose  $a \in X$ and let $x \in I_{P^{-1}a}$ be given i.e. $\alpha_{a}(x)=0$. Let $(P^{-1}b_{n})_n$ be a sequence converging to $X$. Fix $s \in Int(P)$ and write $x=\alpha_{s}(y)$.  Then $P^{-1}b_{n}s \to Xs$. Since $as \in Xs $, it follows that there exists $c_{n} \in P$ such that $c_{n}^{-1}b_{n}s \to as$. Since $as \in Int(P)$, it follows that $d_{n}:=c_{n}^{-1}b_{n}s \in Int(P)$ eventually. Then $b_{n}s=c_{n}d_{n}$.
Now
\begin{align*}
||\alpha_{b_n}(x)||&=||\alpha_{b_{n}s}(y)|| \\
                       & =||\alpha_{c_{n}d_{n}}(y)|| \\
                        & \leq ||\alpha_{d_n}(y)|| \\
                        & \to ||\alpha_{as}(y)|| \\
                        & = ||\alpha_{a}(x)|| \\
                        &\to 0
                                         \end{align*}
Thus $I_{P^{-1}a} \subset I_{X}$ if $a \in X$ or equivalently if $P^{-1}a \subset X$. If $P^{-1}a \subset X$, denote the map $A/I_{P^{-1}a} \ni x+I_{P^{-1}a} \to x+I_{X} \in A/I_{X}$ be denoted by
$\pi_{X,P^{-1}a}$. For $X \in \Omega$,  let $\pi_{X,X}$ be the identity map on $A_{X}$.

\begin{ppsn}
\label{K-theory-surjectivity}
Let $A$ be a $C^{*}$-algebra and $\alpha:P \to End(A)$ be a surjective action. If the semigroup $P$ satisfies condition \textbf{(H)} then $K_{i}(A \rtimes_{red,\alpha} P)=0$ for $i=0,1$.
\end{ppsn}
\textit{Proof.} Again by Theorem \ref{K-theory}, it is enough to prove that $\Gamma_{0}(\Omega_{0},\cla)$ has trivial $K$-theory where $\cla$ is the bundle defined in Section 6. Let $(\phi_{t})_{[0,1]}$ be a homotopy given by Condition \textbf{(H)}.

Let $t \in [0,1]$ be fixed. For $f \in \Gamma(\Omega,\cla)$ and $t \in [0,1]$, let $\pi_{t}(f)$ be the section defined by $\pi_{t}(f)(X):=\pi_{X,\phi_{t}(X)}(f(\phi_{t}(X)))$. We claim that $\pi_{t}(f)$ is continuous if $f \in \Gamma(\Omega,\cla)$. The claim is clear when $f$ is of the form $f(X):=(a+I_{X})\psi(X)$ where $a \in A$ and $\psi \in C(\Omega)$. But then the linear span of functions of the above form are dense in $\Gamma(\Omega,\cla)$ and $\pi_{t}$ is contractive. This proves the claim.

Thus we obtain a homotopoy of $*$-algebra homomorphisms $(\pi_{t})_{t \in [0,1]}: \Gamma(\Omega,\cla) \to \Gamma(\Omega,\cla)$. Now $(1)$ of Condition \textbf{(H)} implies that $\pi_{t}$ leaves $\Gamma_{0}(\Omega_{0},\cla)$ invariant for every $t$. By $(3)$ of Condition \textbf{(H)}, the restrictions of $(\pi_{t})_{t \in [0,1]}$ to $\Gamma_{0}(\Omega_{0},\cla)$ is a homotopy of $*$-algebra homomorphisms connecting the identity map with the zero map. Thus the $K$-theory of $\Gamma_{0}(\Omega_{0},\cla)$ vanishes. This completes the proof. \hfill $\Box$

\begin{rmrk}
The vanishing of $K$-groups when $P:=[0,\infty)$ and $G:=\mathbb{R}$ is proved in complete generality in \cite{KS97}. We remark that $P$ satisfies Condition \textbf{(H)}. In this case, the unit space $\Omega$ of the Wiener-Hopf groupoid can be identified with $[0,\infty]$, the one-point compactification of $[0,\infty)$. The map $[0,\infty] \ni x \to (-\infty,x]$ is a homeomorphism. Moreover the right action of $P$ on $[0,\infty]$ is by translation with the understanding that $\infty + x = \infty $ for $x \in P$.  For $t \in [0,1]$, let $\phi_{t}:[0,\infty] \to [0,\infty]$ be defined by \[\phi_t(x):=\frac{(1-t)x}{\sqrt{(1-(1-t)^{2})x^{2}+1}}.\] Then $(\phi_{t})_{t \in [0,1]}$ is the desired homotopy.
\end{rmrk}

It would be an interesting exercise to determine if the positive semigroup of a Heisenberg group considered in \cite{Nica90} satisfies condition \textbf{(H)}. 

\section{The semigroup of positive matrices}
The aim of this section is to show that the semigroup of positive elements associated to a special Euclidean Jordan algebra satisfies condition \textbf{(H)}. Thus the $K$-theory of the Toeplitz algebras
associated to surjective or injective actions of such semigroups vanishes. We need a few preliminaries regarding the action of the self-adjoint matrices on the unitary group.  The Wiener-Hopf groupoid associated
to such semigroups were described in \cite{Renault_Muhly}. However we obtain a pleasant picture of the unit space in terms of unitary matrices which is helpful for the homotopy considerations.

\subsection{Cayley transform and the M\"{o}bius action}
In what follows, let $\mathcal{H}$ be a finite dimensional complex Hilbert space. Denote the algebra of linear operators on $\clh$ by  $B(\mathcal{H})$.  If $A,B \in B(\mathcal{H})$ are self-adjoint, we write $A \geq B$ if $A-B$ is positive. We write $A>B$ if $A-B$ is  positive and invertible. Also for $A \in B(\clh)$,  let $\sigma(A)$ be its spectrum. Let 
\begin{align*}
\mathcal{S}:&= \{A \in B(\clh): \textrm{ $A$ is self-adjoint}\}, \textrm{~and} \\ 
P:& = \{A \in \mathcal{S}: \textrm{ $A$ is positive} \}.
\end{align*}
Then $\mathcal{S}$ is a real vector subspace of $B(\clh)$ and $P$ is a closed convex cone in $\mathcal{S}$. Also note that $P-P= \mathcal{S}$. It is easily verifiable that 
\[
 Int(P):=\{A \in \mathcal{S}: A >0\}.
\]
Also observe that $Int(P)$ is dense in $P$. Let $\mathcal{U}(\clh)$ be the set of unitaries on $\clh$. The Cayley transform \[\mathcal{S} \ni A \rightarrow \frac{A+i}{A-i} \in \mathcal{U}(\clh)\] is injective and its image is $\{U \in \mathcal{U}(\clh): 1 \notin \sigma(U) \}$ which  is dense  in $\mathcal{U}(\clh)$.

Let $\mathbb{T}_{+}:=\{ z \in \mathbb{T}: Im(z) \geq 0 \}$. It is easy to check that the image of $P$, under the Cayley transform, is $\{ U \in \mathcal{U}(\clh): \sigma(U) \subset \mathbb{T}_{+}\setminus\{1\} \}$. Denote the closure of $P$ in $\mathcal{U}(\clh)$ by $X$. Then \[
                                                                                                                                                                                                                                                                                                                      X= \{ U \in \mathcal{U}(\clh): \sigma(U) \subset \mathbb{T}_{+}\}.
                                                                                                                                                                                                                                                                                                                     \]  The following lemma is useful in proving this and we leave its proof to the reader.

\begin{lmma}
\label{Continuity of the spectrum}
 Let $K$ be a compact subset of $\mathbb{C}$. Supposet $U_{n}$ is a sequence in $\mathcal{U}(\clh)$ such that $U_{n}$ converges to $U$. If $\sigma(U_{n}) \subset K$ then $\sigma(U) \subset K$.
\end{lmma}

Our aim is to show that the action of $\mathcal{S}$ on itself by addition extends to an action on $\mathcal{U}(\clh)$. We start with a lemma.

\begin{lmma}
 Let $U \in \mathcal{U}(\clh)$ and $B \in \mathcal{S}$. Then $BU+2i-B$ is invertible.
\end{lmma}
\textit{Proof.} By taking adjoints, it is enough to show that $UB-2i-B$ is invertible. Suppose $\xi \in \clh$ be such that $(UB-2i-B)\xi=0$. Then $(U-1)B\xi=2i\xi$. Let $\eta=B\xi$ and let $\displaystyle U=\sum_{k=1}^{\ell}\lambda_{i}E_{i}$ be the spectral decomposition of $U$. Since $B$ is self-adjoint, it follows, from the equality $<(U-1)\eta,\eta>=2i<\xi,B \xi>$, that $Re(<(U-1)\eta,\eta>)=0$. Now 
\begin{displaymath}
 Re(<(U-1)\eta,\eta>)= \sum_{k=1}^{\ell}Re(\lambda_{k}-1)<E_{k}\eta,\eta> .
\end{displaymath}
Since $\lambda_{k} \in \mathbb{T}$, it follows that for every $k$, $Re(\lambda_{i}-1) \leq 0$. Also if $\lambda_{i} \neq 1$, $Re(\lambda_{i}-1)<0$. Hence $E_{i}\eta=0$ if $\lambda_{i} \neq 1$. This shows that $(U-1)\eta=0$. Since $(U-1)\eta = 2i \xi$, it follows that $\xi=0$. Thus  $UB-2i-B$ is one-one and hence invertible.  \hfill $\Box$.

For $U \in \mathcal{U}(\clh)$ and $B \in \mathcal{S}$, let \[U \boxplus B:= \big((2i+B)U-B\big)\big(BU+2i-B\big)^{-1}.\] If $A \in \mathcal{S}$, let $U_{A}=(A+i)(A-i)^{-1}$. It is easy to check that  $U_{A} \boxplus B= U_{A+B}$. It then follows that, by the density of $\mathcal{S}$ in $\mathcal{U}(\clh)$, $U \boxplus B \in \mathcal{U}(\clh)$ for $U \in \mathcal{U}(\clh)$ and $B \in \mathcal{S}$. This way, one obtains a right action of $\mathcal{S}$ on $\mathcal{U}(\clh)$.

\begin{rmrk}
 For $B \in \mathcal{S}$, the  action of $B$ on $\mathcal{U}(\clh)$ is via the M\"{o}bius transformation $\begin{pmatrix}
                                                                                                    2i+B & -B \\
                                                                                                     B & 2i-B
                                                                                                   \end{pmatrix}$
\end{rmrk}

\subsection{Special Jordan algebras}
For the rest of the paper, let $V$ be a  real vector subspace of $B(\clh)$ such that
\begin{enumerate}
  \item [(1)] for $a,b \in V$, $a \circ b = \frac{1}{2}(ab+ba) \in V$,
 \item [(2)] the identity $1 \in V$, and
 \item [(3)] If $a \in V$ then $a$ is self-adjoint.
\end{enumerate}
For example, we can let $V$ to be the subspace of self-adjoint operators.

\begin{rmrk}
The multiplication $\circ$ is called a Jordan product on $V$ and the pair $(V,\circ)$ is an example of a special Euclidean Jordan algebra. 
For an introduction to Jordan algebras and for other examples, we refer the reader to \cite{Faraut} and \cite{Mccrimmon}.  However, we do not presume any knowledge of Jordan algebra to understand this section.

\end{rmrk}

Let $Q:= \{a \in V: a \geq 0\}$. Note that the real algebra generated by $\{1,a\} \subset V$ for every $a \in V$. Hence we have functional calculus at our disposal i.e. if $a \in V$ and $f \in C_{0}(\mathbb{R})$ is a real-valued function then $f(a) \in V$. Thus $V=Q-Q$. Also note that $Int(Q)$, in $V$, is $\{a \in Q: a >0\}$. Moreover $V \cong \mathbb{R}^{n}$ for some $n$ and hence $V$ satisfies the Connes-Thom isomorphism property with period $n$.

Embedd $Q$ in $\mathcal{U}(\clh)$ via the Cayley transform and let $Z$ be the closure of $Q$.  Then  $ \{U \in \mathcal{U}(\clh): \sigma(U) \subset \mathbb{T}_{+}\}$ contains $Z$. We endow $Z$ with the subspace topology inherited from the norm topology on $\mathcal{U}(\clh)$. Then $Z$ is compact. 
\begin{lmma}
 \label{functional calculus}

  If $f:\mathbb{T} \to \mathbb{R}$ is continuous and $U \in Z$ then $f(U) \in V$. Similarly if  $g:\mathbb{R} \to \mathbb{T}_{+}$ be continuous and $A \in Q$. Then $g(A) \in Z$.
\end{lmma}
\textit{Proof.} Let $f:\mathbb{T} \to \mathbb{R}$ be continuous. Since $Q$ is dense in $Z$, it is enough to show that   $f(U)\in V$ where $\displaystyle U =\frac{A+i}{A-i}$ for some $A \in Q$. Fix such an $U$.  Define $g:\mathbb{R} \to \mathbb{R}$ by $g(x):=f(\frac{x+i}{x-i})$. Then $f(U)=g(A)$. Since the real algebra generated by $\{1,A\}$ is contained in $V$, it follows that $g(A) \in V$ and thus $f(U) \in V$.

Now let $g:\mathbb{R} \to \mathbb{T}_{+}$ be continous. Consider an element $A \in Q$. Let $A:=\sum_{i=1}^{n}\lambda_{k}E_{k}$ be the spectral decomposition of $A$. Since $E_{k}$ is a polynomial in $A$ with real co-efficients, it follows that $E_{k} \in V$ for every $k$. Now $g(\lambda_{k}) \in \mathbb{T}_{+}$ for every $k$ and $\mathbb{T}_{+}\setminus \{1\}$ is dense in $\mathbb{T}_{+}$. For every $k$, choose a sequence $(\lambda_{k}^{(n)})$ of positive real numbers such that $\displaystyle \frac{\lambda_{k}^{(n)}+i}{\lambda_{k}^{(n)}-i}$ converges to $g(\lambda_{k})$. Set $A^{(n)}:= \sum_{k}\lambda_{k}^{(n)}E_{k}$ and let $U^{(n)}$ be its Cayley transform. Then $A^{(n)} \in Q$ and $U^{(n)}$ converges to $g(A)$. Since $Z$ is closed, it follows that $g(A) \in Z$. This completes the proof. \hfill $\Box$

Since $Q$ is dense in $Z$ and the M\"{o}bius action of $Q$ on $\mathcal{U}(\clh)$ leaves $Q$ invariant, it follows that $Q$ leaves $Z$ invariant i.e. $U \boxplus A \in Z$ for $U \in Z$ and $A \in Q$. Let \[
                                                                                                                                                                                                        Z_{0}:= \{U \boxplus A: U \in Z, A \in Int(Q)\}.
                                                                                                                                                                                                       \]

We claim that the map $Z \times Int(Q) \ni (U,A) \to U \boxplus A \in Z$ is open. By Theorem 4.3 in \cite{Jean_Sundar} , it is enough to show that for every $B \in Q$, $Z_{0} \boxplus B$ is open in $Z$.  Let $B \in Q$ be given. Note that $Z_{0} \boxplus B= \displaystyle \bigcup_{A \in Int(Q)}Z_{0} \boxplus (A +B)$. 
(For, given $C \in Int(Q)$, there exists $C_{1},C_{2} \in Int(Q)$ such that $C=C_{1}+C_{2}$).  Hence it is enough to show that $Z_{0} \boxplus B$ is open whenever $B \in Int(Q)$. 

\begin{rmrk}
 Observe that if $a,b \in \mathbb{R}$ are non-negative then $a > b$ if and only if $Re(\frac{a+i}{a-i}) > Re(\frac{b+i}{b-i})$.
\end{rmrk}

First we need the following lemma.

\begin{lmma}
 We have $Z_{0}=\{ U \in Z: -1 \notin \sigma(U)\}$. Also the set $Z_{0}$ is open in $Z$.
\end{lmma}
\textit{Proof.} Let $U \in Z$ and $-1 \notin \sigma(U)$. Let $U=\sum_{k}\lambda_{k}E_{k}$ be the spectral decomposition. By Lemma \ref{functional calculus}, it follows that $E_{k} \in Q$ for every $k$. For every $k$, $\lambda_{k} \neq -1$. Choose $a>0$ such that $ Re(\lambda_{k})>\frac{a^{2}-1}{a^{2}+1}= Re(\frac{a+i}{a-i})$ for every $k$. For every $k$, choose a sequence $(\lambda_{k}^{(n)}) \in \mathbb{T}_{+} \setminus \{1\}$ such that $Re(\lambda_{k}^{(n)}) > Re(\frac{a+i}{a-i})$, and $(\lambda_{k}^{(n)})$ converges to $\lambda_{k}$. Write $\lambda_{k}^{(n)}=\frac{\mu_{k}^{(n)}+i}{\mu_{k}^{(n)}-i}$. Set $A^{(n)}:= \sum_{k} \mu_{k}^{(n)}E_{k}$ and $U^{(n)}$ be its Cayley transform. Then $A^{(n)} \in Q $ and $A^{(n)}-a \in Q$. Thus $U^{(n)}\boxplus (-a) \in Z$ and hence its limit $V:=U\boxplus (-a) \in Z$.  Thus $U=V \boxplus a$. This proves that  $U \in Z_{0}$.

Let $U \in Z_{0}$. Let $U=\widetilde{U} \boxplus B$ with $\widetilde{U} \in Z$ and $B \in Int(Q)$. Choose a sequence $A_{n} \in Q$ such that its Cayley transform $U_{n}$ converges to $\widetilde{U}$. Let $b$ be the smallest eigen value of $B$. Since $B$ is invertible, $b>0$. Then $\sigma(A_{n}+B) \subset [b,\infty)$. Hence $\sigma(U_{n}\boxplus B) \subset \{z \in \mathbb{T}_{+}: Re(z) \geq Re(\frac{b+i}{b-i})\}$. By Lemma \ref{Continuity of the spectrum}, it follows that $\sigma(\widetilde{U} \boxplus B) \subset \{z \in \mathbb{T}_{+}: Re(z) \geq Re(\frac{b+i}{b-i})\}$. Hence $-1 \notin \sigma(U)$. 

The openness of $Z_{0}$ in $Z$ is now clear. This completes the proof. \hfill $\Box$

We can identify $\mathbb{T}_{+}\setminus \{-1\}$ with $[0,\infty)$ via the map $[0,\infty) \ni x \to \frac{-x+i}{x+i}$. Thus we can identify $Z_{0}$ with $Q$ homeomorphically. Let $\psi:Q \to Z_{0}$ be the map defined by \[\psi(A)=\frac{-A+i}{A+i}.\] Then $\psi$ is a bijection. Its inverse is given by $\psi^{-1}(U)=i(1-U)(1+U)^{-1}$ for $U \in Z_{0}$. Note that $\psi(A) \in Z$ by Lemma \ref{functional calculus}. Thus $\psi$ is well defined. Again by Lemma \ref{functional calculus}, $\psi^{-1}$ is well defined.

By direct verification, it follows that for $A \in Q$, $B \in Q$, \[\psi^{-1}(\psi(A) \boxplus B)=A(BA+1)^{-1}.\]
Let $B \in Int(Q)$. Since $Z_{0}$ is open in $Z$, $Z_{0}\boxplus B$ is open in $Z$ if and only if $Z_{0}\boxplus B$ is open in $Z_{0}$. It is thus enough to show that the range of the map $Q \ni A \to A(BA+1)^{-1} \in Q$ is open in $Q$.

\begin{lmma}
 Let $B \in Int(Q)$ be given. Then the range of the map $Q \ni A \to A(BA+1)^{-1}$ is $\{C \in Q: C < B^{-1}\}$. Moreover $\{C \in Q: C<B^{-1}\}$ is open in $Q$
\end{lmma}
\textit{Proof.} Let $A \in Q$ be given and let $C=A(1+BA)^{-1}$.  Observe that 
\begin{align*}
 B^{-1}-C & = B^{-1}-A(1+BA)^{-1} \\
          & = B^{-1}((1+BA)-BA)(1+BA)^{-1}\\
          & = B^{-1}(1+BA)^{-1}
\end{align*}
Now since $(1+BA)B=B+BAB >0$, it follows that $B^{-1}(1+BA)^{-1} > 0$. Hence $B^{-1}-C>0$.

Now let $C \in Q$ be such that $C<B^{-1}$. Let $D \in Int(Q)$ be such that $B^{-1}-C=D$. Now $1-CB=DB$ and hence $1-CB$ is invertible. Set $A=(1-CB)^{-1}C$. Note that $A$ is positive, for $A=B^{-1}D^{-1}C=(C+D)D^{-1}C=C+CD^{-1}C \geq 0$.

Next we show that $A \in Q$. First we prove that $A \in Q$ if $C$ is invertible. The equality $A=C(1+BA)$, together with the assumption that $C$ is invertible, implies that $C^{-1}=(1+BA)A^{-1}=A^{-1}+B$. Since $C^{-1},B \in Q$, it follows that $A^{-1} \in Q-Q$ and hence $A \in Q-Q$.  We have already shown that $A$ is positive. Thus $A \in Q$. 

Now let $C \in Q$ be such that $C<B^{-1}$ and write $B^{-1}=C+D$. Since $D$ is invertible, $D-\frac{1}{n}>0$ eventually. Set $C_{n}=C+\frac{1}{n}$. Then $C_{n}$ is invertible and converges to $C$. Also note that $B^{-1}-C_{n}= D -\frac{1}{n}>0$. Thus $(1-C_{n}B)^{-1}C_{n} \in Q$. Now the closedness of $Q$ and the convergence of $(1-C_{n}B)^{-1}C_{n}$ to $A:=(1-CB)^{-1}C$ implies that $A \in Q$. The equality $C(1+BA)=A$ implies $C=A(1+BA)^{-1}$.

Now we show that $\{C \in Q: C<B^{-1}\}$ is open in $Q$. Let $C \in Q$ be such that $C<B^{-1}$ and write $B^{-1}=C+D$. Then $D>0$. Thus there exists $\epsilon >0$ such that $<D\xi,\xi> \geq \epsilon$ if $||\xi||=1$. Now the ball of radius $\frac{\epsilon}{2}$ and with centre $C$ is contained in $\{C \in Q: C<B^{-1}\}$. This completes the proof. \hfill $\Box$

For $U \in Z$, let \[Q_{U}:=\{A \in V: U \boxplus A \in Z\}.\] Now we verify the following conditions.
\begin{enumerate}
 \item[(A1)] The $Q$-orbit of $-1 \in Z$ is dense in $Z$,
 \item[(A2)] $Q_{-1}=Q$, and
 \item[(A3)] For $U,V \in Z$, if $Q_{U}=Q_{V}$ then $U=V$.
\end{enumerate}
These conditions are considered in Section 5 of \cite{Jean_Sundar}. (A1) is just the definition of $Z$. To verify these conditions, we use another description of $Z$, the one used in \cite{Renault_Muhly}.

Let $U$ be a unitary on $\clh$. Let $E$ be the projection onto the eigen space $Ker(U-1)$. On $E^{\perp}$, $U-1$ is invertible. Thus $U|_{E^{\perp}}$ comes from a self-adjoint operator on $1-E$, via the Cayley transform. Thus $U$ can be described by a pair $(E,A)$ where $E$ is a projection and $A$ is a self-adjoint operator on $\clh$ such that $(1-E)A(1-E)=A$.

Let \[\widetilde{Z}:=\{(E,A): \textrm{~$E$ is a projection in $Q$, $A \in Q$ such that $(1-E)A(1-E)=A$}\}.\]

For $(E,A) \in \widetilde{Z}$, let $\displaystyle U_{(E,A)}=E+(1-E)\big (\frac{A+i}{A-i} \big)(1-E)$. By functional calculus ( Lemma \ref{functional calculus} ), it follows that $U_{(E,A)} \in Z$. Moreover $\widetilde{Z} \in (E,A) \to U_{(E,A)} \in Z$ is a bijection. Also if $(E,A) \in \widetilde{Z}$ and $B \in V$ then \[U_{(E,A)}\boxplus B=U_{(E,A+ (1-E)B(1-E))}~.\] We leave these verification to the reader. To prove the final assertion, decompose $\clh$ as $E \oplus (1-E)$ and write $B$ as a $2 \times 2$ matrix with respect to this decomposition. 

For $(E,A) \in Z$, we simply write $Q_{(E,A)}$ instead of $Q_{U_{(E,A)}}$. From this it follows that for $(E,A) \in \widetilde{Z}$, \[
                                                                       Q_{(E,A)}= \{ B \in V: A+(1-E)B(1-E) \geq 0 \}.
                                                                      \]
Now $Q_{-1}=Q_{(0,0)}=\{B \in V: B \geq 0\}=Q$. This proves (A2). 

Let $(E,A), (F,B) \in Z$ be such that $Q_{(E,A)}=Q_{(F,B)}$. First observe that $\alpha E \in Q_{(E,A)}$ for every $\alpha \in \mathbb{R}$. Hence $\alpha E \in Q_{(F,B)}$ for $\alpha \in \mathbb{R}$. Thus $B + \alpha (1-F)E(1-F) \geq 0$ for $\alpha \in \mathbb{R}$. This implies that $\xi \in Ker(F)$, $\langle B\xi,\xi \rangle +\alpha \langle  E\xi,\xi \rangle \geq 0$ for every $\alpha \in \mathbb{R}$. Hence $\langle E\xi,\xi \rangle=0$ or in other words $(1-E)\xi=\xi$. Hence $1-F \leq 1-E$ or $E \leq F$. Similarly $F \leq E$. Thus $E=F$. 

Now observe that $-A \in Q_{(E,A)}$. Hence $-A \in Q_{(F,B)}$ which implies that $B-A \geq 0$. Thus $A \leq B$. Similarly $B \leq A$. Hence $A=B$. This proves (A3).

Now we prove our main theorem.  Let $\Omega$ be the unit space of the Wiener-Hopf groupoid associated to the semigroup $Q$. Then by Proposition 5.1 of \cite{Jean_Sundar}, it follows that the map $Z \ni U \to Q_{U}^{-1} \in \Omega$ is a $Q$-equivariant 
homeomorphism. Via this identification, we can identify the boundary $\partial \Omega$, also denoted $\partial Z$,  with $\{U \in Z: -1 \in \sigma(U)\}$.

\begin{rmrk}
\label{order relation}
We identify $Q$ with a subset of $Z$ via the Cayley transform $A \to \frac{A+i}{A-i}$. Let $U \in Z$. Then
\begin{enumerate}
\item[(1)] $U \in Q$ if and only if $1 \notin \sigma(U)$.
\item[(2)] Let $U:=\sum_{i=1}^{k}\lambda_{i}E_{i}$ be the spectral decomposition and let $V:=\sum_{i=1}^{k}\mu_{k}E_{k} \in Z$. Then $V \in Q$ if and only if $\mu_{k} \neq 1$ for every $k$. Moreover if we write $V:=\frac{A+i}{A-i}$ with $A$ self-adjoint then $A \in Q_{U}^{-1}$ if and only if 
$Re(\mu_{k}) \leq Re(\lambda_{k})$ for every $k$. We leave this verification to the reader.
\end{enumerate}

\end{rmrk}

\begin{thm}
The semigroup $Q$ satisfies condition \textbf{(H)}.
\end{thm}
\textit{Proof.} Let $g$ be the inverse of the map $[0,\pi] \ni t \to e^{it} \in \mathbb{T}_{+}$. For $t \in [0,1]$, let $\phi_{t}:Z \to Z$ be defined by $\phi_{t}(U):=e^{i(1-t)g(U)+it\pi}$. Clearly if $U \in \partial Z$ then
$\phi_{t}(U) \in \partial Z$. Moreover $\phi_{0}:=id_{Z}$ and $\phi_{1}(U)= -1 \in \partial Z$ for $U \in Z$. 

We leave it to the reader to verify using Remark \ref{order relation} that for $t \in (0,1]$, $\phi_{t}(U) \in Q$ and  $U \in Z$ and if we write $\phi_{t}(U):=\frac{A_{t,U}+i}{A_{t,U}-i}$ with $A_{t,U}$ self-adjoint then
$A_{t,U} \in Q_{U}^{-1}$. The fact that $A_{t,U} \in Q_{U}^{-1}$ follows from Remark \ref{order relation} and the fact that cosine is decreasing on $[0,\pi]$.  Thus $(\phi_{t})_{t \in [0,1]}$ is the desired homotopy for the semigroup $Q$ to satisfy condition \textbf{(H)}. \hfill $\Box$

\nocite{Nica92} 
 \bibliography{references}
\bibliographystyle{amsalpha}

\noindent
{\sc S. Sundar}
(\texttt{sundarsobers@gmail.com})\\
         {\footnotesize  Chennai Mathematical Institute, H1 Sipcot IT Park, \\
Siruseri, Padur, 603103, Tamilnadu, INDIA.}

\end{document}